\documentclass{article}
\usepackage{amsfonts,amsmath,amssymb}
\usepackage{graphics, epsfig}
\usepackage{color}
\allowdisplaybreaks
\makeatletter
\newcommand{\bpsi}{\boldsymbol\psi}
\newcommand{\AppendicesFromNowOn}
  {\renewcommand \thesection {\@Alph\c@section}  % Changed from book.cls
\newcommand\Appendix{\@startsection {section}{1}{\z@}%
                                   {-3.5ex \@plus -1ex \@minus -.2ex}%
                                   {2.3ex \@plus.2ex}%
                                   {\noindent\normalfont\Large\bfseries Appendix }}
   \setcounter{section}{0}
  }
\makeatother
\begin{document}
%%%%%%%%%%%%%%%%%%%%%%%%%%%%%%%%%%%%%%%%%%%%%%%%%%%
%% if using the Springer package svmno.cls the 
%% theorem envronment is already defined
%% if not remove the comment below
%%%%%%%%%% LEMMA, THM, PROP, SECTION %%%%%%%%%%%%%%%
\newtheorem{theorem}{Theorem}[section]
\newtheorem{proposition}{Proposition}[section]
\newtheorem{lemma}{Lemma}[section]
\newtheorem{corollary}{Corollary}[section]
\newtheorem{remark}{Remark}[section]
\newtheorem{definition}{Definition}[section]
%%\newtheorem{proof}{Proof}
%%%%%%%%%%%%%%%%%%%%%%%%%%%%%%%%%%%%%%%%%%%%%%%%%%%%%%%%%%
\renewcommand{\thesection}{\arabic{section}}
\renewcommand{\theequation}{\thesection.\arabic{equation}}
\renewcommand{\thetheorem}{\thesection.\arabic{theorem}}
\numberwithin{equation}{section}
\numberwithin{theorem}{section}
\numberwithin{proposition}{section}
\numberwithin{lemma}{section}
\numberwithin{remark}{section}
\setcounter{secnumdepth}{3}
%%%%%%%%%%%%%%%%%%%%%%%%%%%%%%%%%%%%%%%%%%%%%%%%%%%%%%%%%%%%%%%
%%%%%%%%%%%%%%%%%%%%%%%%%%%%%%%%%%%%%%%%%%%%%%%%%%%%%%%%%%%%%%%
%%%%%%%%%%%%%%%%%%%%%%%%%%%%%%%%%%%%%%%%%%%%%%%%%%%%%%%%%%%%%%%
%%%% FORMATTING MACROS %%%%%%
\newcommand{\cl}{\centerline}
\newcommand{\sms}{\smallskip}
\newcommand{\ms}{\medskip}
\newcommand{\bs}{\bigskip}
\newcommand{\noi}{\noindent}
\newcommand{\itl}[1]{\textit{#1}}
\newcommand{\blf}[1]{\textbf{#1}}
\newcommand{\dsty}{\displaystyle}
\newcommand{\txty}{\textstyle}
\newcommand{\ssty}{\scriptstyle}
\newcommand{\tty}{\texttt}

%%%%%%%%%% SPECIAL SYMBOLS %%%%%%%%%%%%

\newcommand\Par{\mathhexbox278\,}

%%%% GREEK LETTERS  %%%%%%%%

\newcommand{\al}{\alpha}
\newcommand{\Al}{\Alpha}
\newcommand{\be}{\beta}
\newcommand{\Be}{\Beta}
\newcommand{\Gm}{\Gamma}
\newcommand{\gm}{\gamma}
\newcommand{\dl}{\delta}
\newcommand{\Dl}{\Delta}
\newcommand{\lm}{\lambda}
\newcommand{\Lm}{\Lambda}
\newcommand{\kp}{\kappa}
\newcommand{\varep}{\varepsilon}
\newcommand{\eps}{\epsilon}
\newcommand{\vp}{\varphi}
\newcommand{\sig}{\sigma}
\newcommand{\Sig}{\Sigma}
\newcommand{\om}{\omega}
\newcommand{\Om}{\Omega}
\newcommand{\uom}{\mbox{\boldmath$\omega$}}
\newcommand{\btau}{\mbox{\boldmath$\tau$}}
\newcommand{\bnu}{\mbox{\boldmath$\nu$}}
\newcommand{\up}{\upsilon}
\newcommand{\z}{\zeta}

%%% SPECIAL MATH SYMBOLS %%%%%%%%%%%%

\newcommand{\df}[1]{\buildrel\mbox{\small def}\over{#1}}
\newcommand{\op}[1]{\buildrel\mbox{\tiny o}\over{#1}}
\newcommand{\db}{\prime\prime}
\newcommand{\bsl}{\backslash}
\newcommand{\lnrm}{\|\!|}
\newcommand{\rnrm}{|\!\|}
\newcommand{\lb}{\lbrack\!\lbrack}
\newcommand{\rb}{\rbrack\!\rbrack}
\newcommand\la{\langle}
\newcommand\ra{\rangle}
\newcommand{\ev}{\equiv}
\newcommand{\nev}{\not\equiv}
\newcommand{\nn}{\mathbb{N}}
\newcommand{\qq}{\mathbb{Q}}
\newcommand{\zz}{\mathbb{Z}}
\newcommand{\rr}{\mathbb{R}}
\newcommand{\rn}{\rr^N}
\newcommand{\cc}{\mathbb{C}}
\newcommand{\id}{\mathbb{I}}
\newcommand{\bo}{\mathbb{O}}

\newcommand{\amsb}[1]{\mathbb{#1}}
\newcommand{\mcl}[1]{\mathcal{#1}}
\newcommand{\bl}[1]{\mathbf{#1}}
\newcommand{\ov}[1]{\overline{#1}}
\newcommand{\wt}[1]{\widetilde{#1}}
\newcommand{\wh}[1]{\widehat{#1}}

\newcommand{\llra}{\leftrightarrow}
\newcommand{\lra}{\longrightarrow}
\newcommand{\LLR}{\Longleftrightarrow}
\newcommand{\LRA}{\Longrightarrow}
\newcommand{\LLA}{\Longleftarrow}

%%%% BLACK BOX AND OPEN BOX %%%%%%%%%%%%%%%%

\newcommand{\bbox}{\vrule height.6em width.6em 
depth0em} %%%%% Black Box
\newcommand{\os}{\vbox{\hrule \hbox{\vrule 
height.6em depth0pt 
\hskip.6em \vrule height.6em depth0em}
\hrule}} %%%%%% Open Square

%%%%%%%%%%%%%% OPERATORS %%%%%%%%%%%%%%%%%%%

\newcommand{\Ker}{\operatorname{Ker}}
\newcommand{\Imm}{\operatorname{Im}}
\newcommand{\rank}{\operatorname{rank}}
\newcommand{\dvg}{\operatorname{div}}
\newcommand{\curl}{\operatorname{curl}}
\newcommand{\supp}{\operatorname{supp}}
\newcommand{\essup}{\operatornamewithlimits{ess\,sup}}
\newcommand{\essinf}{\operatornamewithlimits{ess\,inf}}
\newcommand{\essosc}{\operatornamewithlimits{ess\,osc}}
\newcommand{\osc}{\operatornamewithlimits{osc}}
\newcommand{\sign}{\operatorname{sign}}
\newcommand{\loc}{\operatorname{loc}}
\newcommand{\diam}{\operatorname{diam}}
\newcommand{\dist}{\operatorname{dist}}
\newcommand{\card}{\operatorname{card}}
\newcommand{\meas}{\operatorname{meas}}
\newcommand{\spn}{\operatorname{span}}
\newcommand{\dtm}{\operatorname{det}}
%

%%%%%%%% OVER AND UNDER LIMITS %%%%%%%%%%%

\newcommand{\overlim}{\mathop{\overline{\lim}}\limits}
\newcommand{\underlim}{\mathop{\underline{\lim}}\limits}
\newcommand{\ttop}[2]{\genfrac{}{}{0pt}{}{#1}{#2}}
\newcommand{\bcu}{\mathop{\txty{\bigcup}}\limits}
\newcommand{\bca}{\mathop{\txty{\bigcap}}\limits}
\newcommand{\bsu}{\mathop{\txty{\sum}}\limits}
\newcommand{\pro}{\mathop{\txty{\prod}}\limits}

%%%%%%%%%%%  DERIVATIVES %%%%%%%%%%%%%%%%

\newcommand{\pl}{\partial}
\newcommand{\ptt}{\frac{\pl}{\pl t}}
\newcommand{\ppx}{\frac\pl{\pl x}}
\newcommand{\dds}{\frac d{ds}}
\newcommand{\ddt}{\frac d{dt}}

%%%%%%%%%%%%%%%%%%%%%%%%%%%%%%%%%%%%%%%%%%%%%%%%
%%%%%%%%%%%%%  INTEGRALS  %%%%%%%%%%%%%%%
%%%%%%%%%%%%%%%%%%%%%%%%%%%%%%%%%%%%%%%%%%%%%%%%
\newcommand{\intl}{\int\limits}
\newcommand{\iintl}{\iint\limits}
%%%%%%%%%%%%%%%%%%%%%%%%%%%%%%%%%%%%%%%%%%%%%%%%
%%%%%%%%%%%%%  INTEGRAL AVERAGES  %%%%%%%%%%%%%%%
%%%%%%%%%%%%%%%%%%%%%%%%%%%%%%%%%%%%%%%%%%%%%%%%
\def\Xint#1{\mathchoice
    {\XXint\displaystyle\textstyle{#1}}%
    {\XXint\textstyle\scriptstyle{#1}}%
    {\XXint\scriptstyle\scriptscriptstyle{#1}}%
    {\XXint\scriptscriptstyle\scriptscriptstyle{#1}}%
    \!\int}
\def\XXint#1#2#3{\setbox0=\hbox{$#1{#2#3}{\int}$}
    \vcenter{\hbox{$#2#3$}}\kern-0.5\wd0}
\def\bint{\Xint-}
\def\dashint{\Xint{\raise4pt\hbox to7pt{\hrulefill}}}
\def\dashiint{\bint\kern-0.15cm\bint}
% integral averages

\newcommand{\ovl}[3]{\int_{#1}^{#2}\kern-#3pt\raise4pt\hbox to7pt{\hrulefill}\ }

\newcommand{\ovll}[3]{\intl_{#1}^{#2}\kern-#3pt\raise4pt\hbox to7pt{\hrulefill}\ }

\newcommand{\tvl}[2]{\iint_{#1}\kern-#2pt\raise4pt\hbox to7pt{\hrulefill}\ }

%%%%%%%%%%% ANALYSIS MACROS %%%%%%%%%%%%%%%%%%%

%%%% Domain \Om %%

\newcommand{\omt}{\Om_T}
\newcommand{\plo}{\partial\Omega}
\newcommand{\ovo}{\bar{\Om} }

%% C-infinity Local spaces and symbols
%
\newcommand{\ci}[1]{C^\infty\!\left({#1}\right)}
\newcommand{\cio}[1]{C_o^\infty\!\left({#1}\right)}
\newcommand{\lloc}[1]{L_{\loc}\!\left({#1}\right)}
\newcommand{\xy}{|x-y|}

%% Integrals

\newcommand{\intom}{\intl_{\Om}}
\newcommand{\intbo}{\intl_{\plo}}
\newcommand{\inom}{\int_{\Om}}
\newcommand{\inbo}{\int_{\plo}}
\newcommand{\intrn}{\intl_{\rn}}

%%%%%%%%%%%% ENDING THE DOCUMENT %%%%%%%%%%%%%%%%

\newcommand{\bye}{\end{document}}

%%%%%%%%%%%%%%%%%%%%%%%%%%%%%%%%%%%%%%%%%%%%%%%%%

%\hsize=5in
%\vsize=7.5in

%
%% Tag
%
%\def\uptag#1#2$${\iftagsleft@\leqno\else\eqno\fi
%  \hbox{\def\pagebreak{\global\postdisplaypenalty-\@M}%
%  \def\nopagebreak{\global\postdisplaypenalty\@M}\rm(#1\unskip)${}^{#2}$}%
%  $$\postdisplaypenalty\z@\ignorespaces}
%

%\def\btag#1#2$${\iftagsleft@\leqno\else\eqno\fi
%  \hbox{\def\pagebreak{\global\postdisplaypenalty-\@M}%
%  \def\nopagebreak{\global\postdisplaypenalty\@M}\rm(#1\unskip)${}_{#2}$}%
%  $$\postdisplaypenalty\z@\ignorespaces}
%
%\def\ubtag#1#2#3$${\iftagsleft@\leqno\else\eqno\fi
%  \hbox{\def\pagebreak{\global\postdisplaypenalty-\@M}%
%  \def\nopagebreak{\global\postdisplaypenalty\@M}
%  \rm(#1\unskip)${}_{#2}^{#3}$}%
%  $$\postdisplaypenalty\z@\ignorespaces}
%
%\def\ptag#1$${\iftagsleft@\leqno\else\eqno\fi
%  \hbox{\def\pagebreak{\global\postdisplaypenalty-\@M}%
%  \def\nopagebreak{\global\postdisplaypenalty\@M}
%  \rm(#1\unskip)${}^{\prime}$}%
%  $$\postdisplaypenalty\z@\ignorespaces}
%
%% tag wih no brackets
%
%\def\etag#1$${\iftagsleft@\leqno\else\eqno\fi
%  \hbox{\def\pagebreak{\global\postdisplaypenalty-\@M}%
%  \def\nopagebreak{\global\postdisplaypenalty\@M}\rm#1\unskip}%
%  $$\postdisplaypenalty\z@\ignorespaces}
%

%%%%%% DEFINIZIONE DI DATA GIORNO, MESE, ANNO %%%%%%%%%%
%
%\newcommand{\data}{\number\day\space \ifcase\month\or
%Gennaio\or Febbraio\or Marzo\or Aprile\or Maggio\or 
%Giugno\or
%Luglio\or Agosto\or Settembre\or Ottobre\or Novembre\or
%Dicembre\fi,\space\number\year}
%

%%\input dibe_031315
%%%%%%%%%%%%%%%%%%%%%%%%%%%%%%%%%%%%%%%%%%%%%%%%%%%
%\input measure.mac
\input harnack_mono.mac
\DeclareRobustCommand{\rchi}{{\mathpalette\irchi\relax}}
\newcommand{\irchi}[2]{\raisebox{\depth}{$#1\chi$}} % inner command, used by \rchi
\newenvironment{ack}{\medskip{\it Acknowledgement.}}{}
%%%%%%%%%%%%%%%%%%%%%%%%%%%%%%%%%%%%%%%%%%%%%%%%%%%
\title{A Necessary and Sufficient Condition for the 
Continuity of Local Minima of Parabolic Variational 
Integrals with Linear Growth}
%%%%%%%%%%%%%%%%%%%%%%%%%%%%%%%%%%%%%%%%%%%%%%%%%%%
\author{Emmanuele DiBenedetto\footnote{Supported by NSF grant DMS-1265548}\\ 
Department of Mathematics, Vanderbilt University\\  
1326 Stevenson Center, Nashville TN 37240, USA\\
email: {\tt em.diben@vanderbilt.edu}
%%%%%%%%%%%%%%%%%%%%%%%%%%%%%%%%%%%%%%%%%%%%%%%%%%%
\and
%%%%%%%%%%%%%%%%%%%%%%%%%%%%%%%%%%%%%%%%%%%%%%%%%%%
Ugo Gianazza\\
Dipartimento di Matematica ``F. Casorati", Universit\`a di Pavia\\ 
via Ferrata 1, 27100 Pavia, Italy\\
email: {\tt gianazza@imati.cnr.it}
%%%%%%%%%%%%%%%%%%%%%%%%%%%%%%%%%%%%%%%%%%%%%%%%%%%
\and
%%%%%%%%%%%%%%%%%%%%%%%%%%%%%%%%%%%%%%%%%%%%%%%%%%%
Colin Klaus$^*$\\
Department of Mathematics, Vanderbilt University\\  
1326 Stevenson Center, Nashville TN 37240, USA\\
email: {\tt colin.j.klaus@vanderbilt.edu}}
%%%%%%%%%%%%%%%%%%%%%%%%%%%%%%%%%%%%%%%%%%%%%%%%%%%
\date{}
%%%%%%%%%%%%%%%%%%%%%%%%%%%%%%%%%%%%%%%%%%%%%%%%%%%
\maketitle
%%%%%%%%%%%%%%%%%%%%%%%%%%%%%%%%%%%%%%%%%%%%%%%%%%%
%%\vskip.4truecm
%%%%%%%%%%%%%%%%%%%%%%%%%%%%%%%%%%%%%%%%%%%%%%%%%%%
\begin{abstract}
%%%%%%%%%%%%%%%%%%%%%%%%%%%%%%%%%%%%%%%%%%%%%%%%%%%
For proper minimizers of parabolic variational integrals 
with linear growth with respect to $|Du|$, we establish 
a necessary and sufficient condition for $u$ to be continuous 
at a point $\pto$, in terms of a sufficient fast decay of the 
total variation of $u$ about $\pto$ (see (1.4) below).  
These minimizers arise also as {proper} solutions 
to the parabolic $1$-laplacian equation. Hence, the 
continuity condition continues to hold {for such solutions} (\S~\ref{S:3}). 
%%%%%%%%%%%%%%%%%%%%%%%%%%%%%%%%%%%%%%%%%%%%%%%%%%%
\vskip.2truecm
%%%%%%%%%%%%%%%%%%%%%%%%%%%%%%%%%%%%%%%%%%%%%%%%%%%
%%%%%%%%%%%%%%%%%%%%%%%%%%%%%%%%%%%%%%%%%%%%%%%%%%%
\vskip.2truecm
%%%%%%%%%%%%%%%%%%%%%%%%%%%%%%%%%%%%%%%%%%%%%%%%%%%
\end{abstract}
%%%%%%%%%%%%%%%%%%%%%%%%%%%%%%%%%%%%%%%%%%%%%%%%%%%
%%\bigskip
%%%%%%%%%%%%%%%%%%%%%%%%%%%%%%%%%%%%%%%%%%%%%%%%%%%
\section{Introduction}\label{S:1}
%%%%%%%%%%%%%%%%%%%%%%%%%%%%%%%%%%%%%%%%%%%%%%%%%%%
Let $E$ be an open subset of $\rn$, and denote 
by $BV(E)$ the space of functions $v\in L^1(E)$ with 
finite total variation \cite{giusti}
%%%%%%%%%%%%%%%%%%%%%%%%%%%%%%%%%%%%%%%%%%%%%%%%%%%
\begin{equation*}
%%%%%%%%%%%%%%%%%%%%%%%%%%%%%%%%%%%%%%%%%%%%%%%%%%%
\|Dv\|(E):=\sup_{\ttop{\vp\in [C_o^1(E)]^N}{|\vp|\le1}}
\Big\{\langle Dv,\vp\rangle=-\ine v\dvg\vp\,dx\Big\}<\infty.
%%%%%%%%%%%%%%%%%%%%%%%%%%%%%%%%%%%%%%%%%%%%%%%%%%%
\end{equation*}
%%%%%%%%%%%%%%%%%%%%%%%%%%%%%%%%%%%%%%%%%%%%%%%%%%%
Here $Dv=(D_1v,\dots,D_Nv)$ is the vector valued 
Radon measure, representing the distributional 
gradient of $v$. A function $v\in BV_{\loc}(E)$ if
$v\in BV(E^\prime)$ for all open sets $E^\prime\subsetneq E$. 
For $T>0$, let $E_T=E\times(0,T)$, 
and denote  by $\dsty L^1(0,T;BV(E))$ the collection 
of all maps $v:[0,T]\to BV(E)$ such that 
%%%%%%%%%%%%%%%%%%%%%%%%%%%%%%%%%%%%%%%%%%%%%%%%%%%
\begin{equation*}
%%%%%%%%%%%%%%%%%%%%%%%%%%%%%%%%%%%%%%%%%%%%%%%%%%%
v\in L^1(E_T),\qquad  \|Dv(t)\|(E)\in L^1(0,T),
%%%%%%%%%%%%%%%%%%%%%%%%%%%%%%%%%%%%%%%%%%%%%%%%%%%
\end{equation*}
%%%%%%%%%%%%%%%%%%%%%%%%%%%%%%%%%%%%%%%%%%%%%%%%%%%
and the maps
%%%%%%%%%%%%%%%%%%%%%%%%%%%%%%%%%%%%%%%%%%%%%%%%%%%
\begin{equation*}
%%%%%%%%%%%%%%%%%%%%%%%%%%%%%%%%%%%%%%%%%%%%%%%%%%%
(0,T)\ni t\to\langle Dv(t),\vp\rangle
%%%%%%%%%%%%%%%%%%%%%%%%%%%%%%%%%%%%%%%%%%%%%%%%%%%
\end{equation*}
%%%%%%%%%%%%%%%%%%%%%%%%%%%%%%%%%%%%%%%%%%%%%%%%%%%
are measurable with respect to the Lebesgue measure 
in $\rr$, for all $\vp\in [C_o^1(E)]^N$. 
%%%%%%%%%%%%%%%%%%%%%%%%%%%%%%%%%%%%%%%%%%%%%%%%%%%

A function $u\in L^1_{\loc}\big(0,T;BV_{\loc}(E)\big)$ is a local 
parabolic minimizer of the total variation flow 
in $E_T$, if %%there exists $Q\ge1$ such that 
%%%%%%%%%%%%%%%%%%%%%%%%%%%%%%%%%%%%%%%%%%%%%%%%%%%
\begin{equation}\label{Eq:1:1}
%%%%%%%%%%%%%%%%%%%%%%%%%%%%%%%%%%%%%%%%%%%%%%%%%%%
\int_0^T\Big[\ine -u\vp_tdx+\|Du(t)\|(E)\Big]dt\le
\int_0^T \|D(u+\vp)(t)\|(E)dt
%%%%%%%%%%%%%%%%%%%%%%%%%%%%%%%%%%%%%%%%%%%%%%%%%%%
\end{equation}
%%%%%%%%%%%%%%%%%%%%%%%%%%%%%%%%%%%%%%%%%%%%%%%%%%%
for all non-negative $\vp\in C_o^\infty(E_T)$. The notion has been 
introduced in \cite{BDM}. It is a parabolic 
version of the elliptic local minima of total variation flow 
as introduced in \cite{HK}.
%%%%%%%%%%%%%%%%%%%%%%%%%%%%%%%%%%%%%%%%%%%%%%%%%%%%%%
\subsection{The Main Result}\label{S:1:1}
%%%%%%%%%%%%%%%%%%%%%%%%%%%%%%%%%%%%%%%%%%%%%%%%
Let $B_\rho(x_o)$ denote the ball of radius $\rho$ 
about $x_o$. If $x_o=0$, write $B_\rho(x_o)=B_\rho$.
Introduce the cylinders
%%%%%%%%%%%%%%%%%%%%%%%%%%%%%%%%%%%%%%%%%%%%%%%%%%%
%%\begin{equation*}
%%%%%%%%%%%%%%%%%%%%%%%%%%%%%%%%%%%%%%%%%%%%%%%%%%%
$Q_\rho(\theta)=B_\rho\times(-\theta \rho,0]$,
%%%%%%%%%%%%%%%%%%%%%%%%%%%%%%%%%%%%%%%%%%%%%%%%%%%
%%\end{equation*}
%%%%%%%%%%%%%%%%%%%%%%%%%%%%%%%%%%%%%%%%%%%%%%%%%%%
where $\theta$ is a positive parameter to be chosen 
as needed. If $\theta=1$ we write $Q_\rho(1)=Q_\rho$. 
For a point $\pto\in{\bf \rr}^{N+1}$ we let 
$[\pto+Q_\rho(\theta)]$ be the cylinder 
of ``vertex'' at $\pto$ and congruent to 
$Q_\rho(\theta)$, i.e.,
%%%%%%%%%%%%%%%%%%%%%%%%%%%%%%%%%%%%%%%%%%%%%%%%%%%
\begin{equation*}
%%%%%%%%%%%%%%%%%%%%%%%%%%%%%%%%%%%%%%%%%%%%%%%%%%%
[\pto+Q_\rho(\theta)]=B_\rho(x_o)\times(t_o-\theta\rho,t_o],
%%%%%%%%%%%%%%%%%%%%%%%%%%%%%%%%%%%%%%%%%%%%%%%%%%%
\end{equation*}
and we let $\rho>0$ be so small
that $[\pto+Q_\rho(\theta)]\subset E_T$. 
%%%%%%%%%%%%%%%%%%%%%%%%%%%%%%%%%%%%%%%%%%%%%%%%%%
%%%%%%%%%%%%%%%%%%%%%%%%%%%%%%%%%%%%%%%%%%%%%%%%%%%%%%%
\begin{theorem}\label{Thm:1:1} 
%%%%%%%%%%%%%%%%%%%%%%%%%%%%%%%%%%%%%%%%%%%%%%%%%%%%%%%
Let $u\in L^1_{\loc}\big(0,T;BV_{\loc}(E)\big)$ be a 
local parabolic minimizer of the 
total variation flow in $E_T$, satisfying in addition
%%%%%%%%%%%%%%%%%%%%%%%%%%%%%%%%%%%%%%%%%%%%%%%%%%%
\begin{equation}\label{Eq:1:2}
%%%%%%%%%%%%%%%%%%%%%%%%%%%%%%%%%%%%%%%%%%%%%%%%%%%
u\in L^\infty_{\loc}(E_T)\quad\text{ and }
\quad u_t\in L^1_{\loc}(E_T). 
%%%%%%%%%%%%%%%%%%%%%%%%%%%%%%%%%%%%%%%%%%%%%%%%%%%
\end{equation}
%%%%%%%%%%%%%%%%%%%%%%%%%%%%%%%%%%%%%%%%%%%%%%%%%%%
Then, $u$ is continuous at some $\pto\in E_T$, if 
and only if 
%%%%%%%%%%%%%%%%%%%%%%%%%%%%%%%%%%%%%%%%%%%%%%%%%%%
\begin{equation}\label{Eq:1:3}
%%%%%%%%%%%%%%%%%%%%%%%%%%%%%%%%%%%%%%%%%%%%%%%%%%%
\limsup_{\rho\searrow0} {\frac{\rho}{|Q_\rho|}}
\int_{t_o-\rho}^{t_o}\|Du(\cdot,t)\|(B_\rho(x_o))dt=0.
%%%%%%%%%%%%%%%%%%%%%%%%%%%%%%%%%%%%%%%%%%%%%%%%%%%
\end{equation}
%%%%%%%%%%%%%%%%%%%%%%%%%%%%%%%%%%%%%%%%%%%%%%%%%%%
\end{theorem}
%%%%%%%%%%%%%%%%%%%%%%%%%%%%%%%%%%%%%%%%%%%%%%%%%%%
For stationary, elliptic minimizers, condition (\ref{Eq:1:3}) 
has been introduced in \cite{HK}. The stationary version 
of (\ref{Eq:1:3}) implies that $u$ is quasi-continuous 
at $x_o$. For time-dependent minimizers, however, 
(\ref{Eq:1:3}) gives no information on the possible 
quasi-continuity of $u$ at $\pto$.  Condition (\ref{Eq:1:3}), 
is only a measure-theoretical restriction on the speed at which 
a possible discontinuity may develop at $\pto$. 
For this reason our proof is entirely different 
than \cite{HK}, being based instead on a DeGiorgi-type 
iteration technique that exploits precisely such 
a measure-theoretical information.
%%%%%%%%%%%%%%%%%%%%%%%%%%%%%%%%%%%%%%%%%%%%%%%%%%%%%%%%%%%%%%%%
\section{Comments on Boundedness and Continuity}\label{S:2}
%%%%%%%%%%%%%%%%%%%%%%%%%%%%%%%%%%%%%%%%%%%%%%%%%%%%%%%%%%%%%%%%
The theorem requires that $u$ {is} locally bounded and that 
$u_t\in L^1_{\loc}(E_T)$. In the elliptic 
case, local minimizers of the total gradient flow in $E$, are locally 
bounded (\cite[\S~2]{HK}).  This is not the case, in general, 
for parabolic minimizers in $E_T$, even if 
$u_t\in C^\infty_{\loc}\big(0,T;L^1_{\loc}(E)\big)$. Consider the function 
%%%%%%%%%%%%%%%%%%%%%%%%%%%%%%%%%%%%%%%%%%%%%%%%%%%%%%%%%%%%%%%%%%%%%%%%
\begin{equation*}%%\label{Eq:1:5}
%%%%%%%%%%%%%%%%%%%%%%%%%%%%%%%%%%%%%%%%%%%%%%%%%%%
B_1\times(-\infty,1)\ni {(x,t)}\to F(|x|,t)=(1-t)\frac{N-1}{|x|},
\quad\text{ for }\> N\ge3. 
%%%%%%%%%%%%%%%%%%%%%%%%%%%%%%%%%%%%%%%%%%%%%%%%%%%
\end{equation*}
%%%%%%%%%%%%%%%%%%%%%%%%%%%%%%%%%%%%%%%%%%%%%%%%%%%
Denote by $D_aF$ that component of the measure $DF$ 
which is absolutely continuous with respect to the 
Lebesgue measure in $\rn$.  One verifies that $DF=D_aF$ and 
$\|DF(t)\|(B_1)=\|D_aF(t)\|_{1,B_1}$.  By direct computation
%%%%%%%%%%%%%%%%%%%%%%%%%%%%%%%%%%%%%%%%%%%%%%%%%%%
\begin{equation*}%%\label{Eq:1:5}
%%%%%%%%%%%%%%%%%%%%%%%%%%%%%%%%%%%%%%%%%%%%%%%%%%%
\int_0^T\int_{B_1}\Big(-F\vp_t+\frac{D_aF}{|D_aF|}\cdot D\vp\Big)dxdt=0,
%%%%%%%%%%%%%%%%%%%%%%%%%%%%%%%%%%%%%%%%%%%%%%%%%%%
\end{equation*}
%%%%%%%%%%%%%%%%%%%%%%%%%%%%%%%%%%%%%%%%%%%%%%%%%%%
for all $\vp\in C_o^\infty\big(B_1\times(0,T)\big)$, {$0<T<1$}. From this
%%%%%%%%%%%%%%%%%%%%%%%%%%%%%%%%%%%%%%%%%%%%%%%%%%%
\begin{equation*}%%\label{Eq:1:5}
%%%%%%%%%%%%%%%%%%%%%%%%%%%%%%%%%%%%%%%%%%%%%%%%%%%
\int_0^T\int_{B_1}\Big(-F\vp_t+\frac{D_aF}{|D_aF|}\cdot D_aF\Big)dxdt
=\int_0^T\int_{B_1}\frac{D_aF}{|D_aF|}\cdot D_a(F-\vp)dxdt,
%%%%%%%%%%%%%%%%%%%%%%%%%%%%%%%%%%%%%%%%%%%%%%%%%%%
\end{equation*}
%%%%%%%%%%%%%%%%%%%%%%%%%%%%%%%%%%%%%%%%%%%%%%%%%%%
which yields
%%%%%%%%%%%%%%%%%%%%%%%%%%%%%%%%%%%%%%%%%%%%%%%%%%%
\begin{equation*}%%\label{Eq:F:1}
%%%%%%%%%%%%%%%%%%%%%%%%%%%%%%%%%%%%%%%%%%%%%%%%%%%
\int_0^T\int_{B_1}\Big(-F\vp_t+|D_aF|\Big)dxdt
\le\int_0^T\int_{B_1}|D_a(F-\vp)|dxdt.
%%%%%%%%%%%%%%%%%%%%%%%%%%%%%%%%%%%%%%%%%%%%%%%%%%%
\end{equation*}
%%%%%%%%%%%%%%%%%%%%%%%%%%%%%%%%%%%%%%%%%%%%%%%%%%%
Thus $F$ is a local, unbounded, parabolic minimizer 
 of the total variation flow. 
%%%%%%%%%%%%%%%%%%%%%%%%%%%%%%%%%%%%%%%%%%%%%%%%%%%
The requirement $u\in L^\infty_{\loc}(E_T)$ could be replaced 
by asking that $u\in L^r_{\loc}(E_T)$ for some $r>N$.  
A discussion on this issue is provided in Appendix~\ref{App:B}.
%%%%%%%%%%%%%%%%%%%%%%%%%%%%%%%%%%%%%%%%%%%%%%%%%
\subsection{On the Modulus of Continuity}\label{S:2:1}
%%%%%%%%%%%%%%%%%%%%%%%%%%%%%%%%%%%%%%%%%%%%%%%%%
While Theorem~\ref{Thm:1:1} gives a necessary and 
sufficient condition for continuity at a given point, 
it provides no information on the modulus of continuity 
of $u$ at $\pto$. Consider the two time-independent functions 
in $B_{\rho}\times(0,\infty)$, for some $\rho<1$: 
%%%%%%%%%%%%%%%%%%%%%%%%%%%%%%%%%%%%%%%%%%%%%%%%%%%%%%
\begin{align*}
%%%%%%%%%%%%%%%%%%%%%%%%%%%%%%%%%%%%%%%%%%%%%%%%%%%%%%
u_1(x_1,x_2)&=
\left\{
%%%%%%%%%%%%%%%%%%%%%%%%%%%%%%%%%%%%%%%%%%%%%%%%%%%%%%
\begin{array}{cl}
%%%%%%%%%%%%%%%%%%%%%%%%%%%%%%%%%%%%%%%%%%%%%%%%%%%%%%
{\dsty \frac1{\ln x_1}}\quad&\text{for}\> x_1>0;\\
{}\\
{\dsty 0}\quad&\text{for}\> x_1=0;\\
{}\\
{\dsty -\frac1{\ln (-x_1)}}\quad&\text{for}\> x_1<0.
%%%%%%%%%%%%%%%%%%%%%%%%%%%%%%%%%%%%%%%%%%%%%%%%%%%%%%
\end{array}\right.
%%%%%%%%%%%%%%%%%%%%%%%%%%%%%%%%%%%%%%%%%%%%%%%%%%%%%%
{}\\
%%%%%%%%%%%%%%%%%%%%%%%%%%%%%%%%%%%%%%%%%%%%%%%%%%%%%%
{}\\
%%%%%%%%%%%%%%%%%%%%%%%%%%%%%%%%%%%%%%%%%%%%%%%%%%%%%%
u_2(x_1,x_2)&=
%%%%%%%%%%%%%%%%%%%%%%%%%%%%%%%%%%%%%%%%%%%%%%%%%%%%%%
\left\{
%%%%%%%%%%%%%%%%%%%%%%%%%%%%%%%%%%%%%%%%%%%%%%%%%%%%%%
\begin{array}{cl}
%%%%%%%%%%%%%%%%%%%%%%%%%%%%%%%%%%%%%%%%%%%%%%%%%%%%%%
\ \sqrt{x_1}\quad&\text{for}\> x_1>0;\\
{}\\
-\sqrt{-x_1}\quad&\text{for}\> x_1\le0.
%%%%%%%%%%%%%%%%%%%%%%%%%%%%%%%%%%%%%%%%%%%%%%%%%%%%%%
\end{array}
%%%%%%%%%%%%%%%%%%%%%%%%%%%%%%%%%%%%%%%%%%%%%%%%%%%%%%
\right.
%%%%%%%%%%%%%%%%%%%%%%%%%%%%%%%%%%%%%%%%%%%%%%%%%%%%%%
\end{align*} 
%%%%%%%%%%%%%%%%%%%%%%%%%%%%%%%%%%%%%%%%%%%%%%%%%%%%%%
Both are stationary parabolic minimizers of the total 
{variation} flow in the sense of (\ref{Eq:1:1})--(\ref{Eq:1:2}), 
over $B_{\frac12}\times(0,\infty)$. We establish this 
for $u_1$, the analogous statement for $u_2$ being analogous.  
Since $u_1\in W^{1,1}(B_{\rho})$, and is time-independent, 
one also has $u\in L^1\big(0,T;BV(B_\rho)\big)$. To verify 
(\ref{Eq:1:1}), one needs to show that
%%%%%%%%%%%%%%%%%%%%%%%%%%%%%%%%%%%%%%%%%%%%%%%%%%%%%%
\begin{equation*}\label{Eq:1:1:elliptic}
%%%%%%%%%%%%%%%%%%%%%%%%%%%%%%%%%%%%%%%%%%%%%%%%%%%%%%
\|Du_1\|(B_{\rho})\le\frac1T\int_0^T \|D(u_1+\vp)
(\cdot,t)\|(B_{\rho})dt\tag{*}
%%%%%%%%%%%%%%%%%%%%%%%%%%%%%%%%%%%%%%%%%%%%%%%%%%%%%%
\end{equation*}
%%%%%%%%%%%%%%%%%%%%%%%%%%%%%%%%%%%%%%%%%%%%%%%%%%%%%%
for all $T>0$, and all $\vp\in C_o^\infty(B_{\rho}\times(0,T))$. 
Let $\mcl{H}^k(A)$ denote the $k$-dimensional Hausdorff 
measure of a Borel set $A\subset \rn$. One checks that 
$\mcl{H}^N([Du_1=0])=0$ and there exists a closed set $K\subset B_\rho$, 
such that $\mcl{H}^{N-1}(K)=0$ and
%%%%%%%%%%%%%%%%%%%%%%%%%%%%%%%%%%%%%%%%%%%%%%%%%%%%%%
\begin{equation*}%%\label{Eq:1:1:elliptic}
%%%%%%%%%%%%%%%%%%%%%%%%%%%%%%%%%%%%%%%%%%%%%%%%%%%%%%
\int_{B_{\rho}-K}\frac{Du_1}{|Du_1|}\cdot D\vp\, dx=0,\quad\text{ for all }\> 
\vp\in C_o^\infty(B_\rho-K). 
%%%%%%%%%%%%%%%%%%%%%%%%%%%%%%%%%%%%%%%%%%%%%%%%%%%%%%
\end{equation*}
%%%%%%%%%%%%%%%%%%%%%%%%%%%%%%%%%%%%%%%%%%%%%%%%%%%%%%
From this, {by Lemma~4 of \cite[\S~8]{bombieri}}, for 
all $\psi\in C_o^\infty(B_{\rho})$, one has 
%%%%%%%%%%%%%%%%%%%%%%%%%%%%%%%%%%%%%%%%%%%%%%%%%%%%%%
\begin{equation*}
%%%%%%%%%%%%%%%%%%%%%%%%%%%%%%%%%%%%%%%%%%%%%%%%%%%%%%
\|Du_1\|(B_{\rho})\le\|D(u_1+\psi)\|(B_{\rho}),
%%%%%%%%%%%%%%%%%%%%%%%%%%%%%%%%%%%%%%%%%%%%%%%%%%%%%%
\end{equation*}
%%%%%%%%%%%%%%%%%%%%%%%%%%%%%%%%%%%%%%%%%%%%%%%%%%%%%%
which, in turn, yields (*).  
The two functions $u_1$ and $u_2$ can be regarded as 
equibounded near the origin. They both satisfy (\ref{Eq:1:3}), 
and exhibit quite different moduli of continuity at the origin.  
This occurrence is in line with a remark of Evans 
(\cite{evans}). A sufficiently smooth minimizer of the 
elliptic functional $\|Du\|(E)$ is a function whose level 
sets are surfaces of zero mean curvature. Thus, if $u$ 
is  a minimizer, so is $\vp(u)$ for all continuous 
monotone functions $\vp(\cdot)$. This implies 
that a modulus of continuity cannot be identified 
solely in terms of an upper bound of $u$. 
%%%%%%%%%%%%%%%%%%%%%%%%%%%%%%%%%%%%%%%%%%%%%%%%%%
\section{Singular Parabolic DeGiorgi Classes}\label{S:3}
%%%%%%%%%%%%%%%%%%%%%%%%%%%%%%%%%%%%%%%%%%%%%%%%%%
%Introduce the cylinders
%%%%%%%%%%%%%%%%%%%%%%%%%%%%%%%%%%%%%%%%%%%%%%%%%%%%
%%%\begin{equation*}
%%%%%%%%%%%%%%%%%%%%%%%%%%%%%%%%%%%%%%%%%%%%%%%%%%%%
%$Q_\rho(\theta)=B_\rho\times(-\theta \rho,0]$,
%%%%%%%%%%%%%%%%%%%%%%%%%%%%%%%%%%%%%%%%%%%%%%%%%%%%
%%%\end{equation*}
%%%%%%%%%%%%%%%%%%%%%%%%%%%%%%%%%%%%%%%%%%%%%%%%%%%%
%where $\theta$ is a positive parameter to be chosen 
%as needed. If $\theta=1$ we write $Q_\rho(1)=Q_\rho$. 
%For a point $\pto\in{\bf \rr}^{N+1}$ we let 
%$[\pto+Q_\rho(\theta)]$ be the cylinder 
%of ``vertex'' at $\pto$ and congruent to 
%$Q_\rho(\theta)$, i.e.,
%%%%%%%%%%%%%%%%%%%%%%%%%%%%%%%%%%%%%%%%%%%%%%%%%%%%
%\begin{equation*}
%%%%%%%%%%%%%%%%%%%%%%%%%%%%%%%%%%%%%%%%%%%%%%%%%%%%
%[\pto+Q_\rho(\theta)]=B_\rho(x_o)\times(t_o-\theta\rho,t_o].
%%%%%%%%%%%%%%%%%%%%%%%%%%%%%%%%%%%%%%%%%%%%%%%%%%%%
%\end{equation*}
%%%%%%%%%%%%%%%%%%%%%%%%%%%%%%%%%%%%%%%%%%%%%%%%%%%
Let $\mcl{C}\big(Q_\rho(\theta)\big)$ 
denote the class of all non-negative, 
piecewise smooth, cutoff functions 
$\z$ defined in $Q_\rho(\theta)$, vanishing 
outside $B_\rho$, such that $\z_t\ge0$ and satisfying 
%%%%%%%%%%%%%%%%%%%%%%%%%%%%%%%%%%%%%%%%%%%%%%%%%%%
\begin{equation*}
%%%%%%%%%%%%%%%%%%%%%%%%%%%%%%%%%%%%%%%%%%%%%%%%%%%
|D\z|+\z_t\in L^\infty\big(Q_\rho(\theta)\big).
%%%%%%%%%%%%%%%%%%%%%%%%%%%%%%%%%%%%%%%%%%%%%%%%%%%
\end{equation*}
%%%%%%%%%%%%%%%%%%%%%%%%%%%%%%%%%%%%%%%%%%%%%%%%%%%
For a measurable function $u:E_T\to\rr$ and $k\in\rr$ set
%%%%%%%%%%%%%%%%%%%%%%%%%%%%%%%%%%%%%%%%%%%%%%%%%%%
\begin{equation*}
%%%%%%%%%%%%%%%%%%%%%%%%%%%%%%%%%%%%%%%%%%%%%%%%%%%%
\ukpm=\{\pm(u-k)\wedge0\}.
%%%%%%%%%%%%%%%%%%%%%%%%%%%%%%%%%%%%%%%%%%%%%%%%%%%
\end{equation*}
%%%%%%%%%%%%%%%%%%%%%%%%%%%%%%%%%%%%%%%%%%%%%%%%%%

The singular, parabolic DeGiorgi class  $[DG]^\pm(E_T;\gm)$ 
is the collection of all measurable maps 
%%%%%%%%%%%%%%%%%%%%%%%%%%%%%%%%%%%%%%%%%%%%%%%%%%%
\begin{equation}\label{Eq:3:1}
%%%%%%%%%%%%%%%%%%%%%%%%%%%%%%%%%%%%%%%%%%%%%%%%%%%%
u\in C_{\loc}\big((0,T);L^2_{\loc}(E)\big)\cap 
L^1_{\loc}\big(0,T;BV_{\loc}(E)\big),
%%%%%%%%%%%%%%%%%%%%%%%%%%%%%%%%%%%%%%%%%%%%%%%%%%%
\end{equation}
%%%%%%%%%%%%%%%%%%%%%%%%%%%%%%%%%%%%%%%%%%%%%%%%%%
satisfying 
%%%%%%%%%%%%%%%%%%%%%%%%%%%%%%%%%%%%%%%%%%%%%%%%%%
\begin{equation}\label{Eq:3:2}
%%%%%%%%%%%%%%%%%%%%%%%%%%%%%%%%%%%%%%%%%%%%%%%%%%
\begin{aligned}
%%%%%%%%%%%%%%%%%%%%%%%%%%%%%%%%%%%%%%%%%%%%%%%%%%
&\sup_{t_o-\theta \rho\le t\le t_o}\int_{B_\rho(x_o)}\ukpm^2\z(x,t)dx\\
&\kern2.0cm +\int_{t_o-\theta \rho}^{t_o}
\|D(\ukpm\z)(\tau)\|(B_\rho(x_o))dt\\
&\le\gm {\iint}_{[\pto+Q_\rho(\theta)]}\big[\ukpm|D\z|
+\ukpm^2|\z_t|\big]dxdt+\\
&\kern2.0cm+\int_{B_\rho(x_o)}
\ukpm^2\z(x,t_o-\theta \rho)dx
%%%%%%%%%%%%%%%%%%%%%%%%%%%%%%%%%%%%%%%%%%%%%%%%%%
\end{aligned}
%%%%%%%%%%%%%%%%%%%%%%%%%%%%%%%%%%%%%%%%%%%%%%%%%%
\end{equation}
%%%%%%%%%%%%%%%%%%%%%%%%%%%%%%%%%%%%%%%%%%%%%%%%%%
for all $[\pto+Q_{\rho}(\theta)]\subset E_T$, 
all $k\in\rr$, and all $\z\in{\cal C}([\pto+Q_\rho(\theta)])$, 
for a given positive constant $\gm$.  The singular DeGiorgi 
classes $[DG](E_T;\gm)$ are defined as 
$[DG](E_T;\gm)=[DG]^+(E_T;\gm)\cap [DG]^-(E_T;\gm)$.
%%%%%%%%%%%%%%%%%%%%%%%%%%%%%%%%%%%%%%%%%%%%%%%%%%
\subsection{The Main Result}\label{S:1:3}
%%%%%%%%%%%%%%%%%%%%%%%%%%%%%%%%%%%%%%%%%%%%%%%%%%
The main result of this note is that the necessary and 
sufficient condition of Theorem~\ref{Thm:1:1} holds for 
functions $u\in DG(E_T;\gm)\cap L^\infty_{\loc}(E_T)$. 
Indeed, the proof of Theorem~\ref{Thm:1:1}, only uses 
the local integral inequalities (\ref{Eq:3:2}). In particular, 
the second of (\ref{Eq:1:2}) is not needed.
%%%%%%%%%%%%%%%%%%%%%%%%%%%%%%%%%%%%%%%%%%%%%%%%%%
\begin{proposition}\label{Prop:3:1}
%%%%%%%%%%%%%%%%%%%%%%%%%%%%%%%%%%%%%%%%%%%%%%%%%%
Let $u$ in the functional classes (\ref{Eq:3:1}), 
be a parabolic minimizer of the total variation flow 
in $E_T$, in the sense of (\ref{Eq:1:1}), satisfying in 
addition (\ref{Eq:1:2}). Then $u\in DG(E_T;2)$. 
%%%%%%%%%%%%%%%%%%%%%%%%%%%%%%%%%%%%%%%%%%%%%%%%%%
\end{proposition}
%%%%%%%%%%%%%%%%%%%%%%%%%%%%%%%%%%%%%%%%%%%%%%%%%%%%%
The proof will be given in Appendix~\ref{App:A}. 
%%%%%%%%%%%%%%%%%%%%%%%%%%%%%%%%%%%%%%%%%%%%%%%%%%%%%
\begin{remark}\label{Rmk:3:1} {\normalfont
%%%%%%%%%%%%%%%%%%%%%%%%%%%%%%%%%%%%%%%%%%%%%%%%%%%%%
Note that in the context of $DG(E_T)$ classes, the 
characteristic condition (\ref{Eq:1:3}), holds with 
no further requirement that $u_t\in L^1_{\loc}(E_T)$. 
The latter however is needed to cast a parabolic 
minimizer of the total variation flow into a 
$DG(E_T)$-class as stated by Proposition~\ref{Prop:3:1}.
%%%%%%%%%%%%%%%%%%%%%%%%%%%%%%%%%%%%%%%%%%%%%%%%%%%%%
}%%
%%%%%%%%%%%%%%%%%%%%%%%%%%%%%%%%%%%%%%%%%%%%%%%%%%%%%
\end{remark}%%
%%%%%%%%%%%%%%%%%%%%%%%%%%%%%%%%%%%%%%%%%%%%%%%%%%%%%
\section{A Singular Diffusion Equation}\label{S:4}
%%%%%%%%%%%%%%%%%%%%%%%%%%%%%%%%%%%%%%%%%%%%%%%%%%%%%
Consider formally, the parabolic $1$-Laplacian equation
%%%%%%%%%%%%%%%%%%%%%%%%%%%%%%%%%%%%%%%%%%%%%%%%%%%
\begin{equation}\label{Eq:4:1}
%%%%%%%%%%%%%%%%%%%%%%%%%%%%%%%%%%%%%%%%%%%%%%%%%%%
u_t-\dvg\Big({\dsty\frac{D u}{|D u|}}\Big)=0
\quad\text{ formally in }\> E_T.
%%%%%%%%%%%%%%%%%%%%%%%%%%%%%%%%%%%%%%%%%%%%%%%%%%%
\end{equation}
%%%%%%%%%%%%%%%%%%%%%%%%%%%%%%%%%%%%%%%%%%%%%%%%%%%
Let $\mcl{P}$ be the class of all Lipschitz continuous, 
non-decreasing functions $p(\cdot)$ defined in $\rr$, with 
$p^\prime$ compactly supported. {Denote by $\mcl{C}(E_T)$} 
the class of all non-negative functions $\z$ defined 
in $E_T$, such that $\z(\cdot,t)\in C_o^1(E)$ for all 
$t\in (0,T)$, and $0\le \z_t<\infty$ in $E_T$.  A function 
$u\in C_{\loc}\big (0,T; L^1(E)\big)$ is a local 
solution to (\ref{Eq:4:1}) if 
%%%%%%%%%%%%%%%%%%%%%%%%%%%%%%%%%%%%%%%%%%%%%%%%%%%
\begin{description} %%%
%%%%%%%%%%%%%%%%%%%%%%%%%%%%%%%%%%%%%%%%%%%%%%%%%%%
\item{\bf a.\ } $p(u)\in L^1_{\loc}\big(0,T;BV(E)\big)$, 
for all $p\in\mcl{P}$;
%%%%%%%%%%%%%%%%%%%%%%%%%%%%%%%%%%%%%%%%%%%%%%%%%%%
\item{\bf b.\ } there exists a vector valued function 
${\bf z}\in [L^\infty (E_T)]^N$ with $\|{\bf z}\|_{\infty, E}\le 1$, 
such that $u_t=\dvg{\bf z}$ in $\mcl{D}^{\prime}(E_T)$;
%%%%%%%%%%%%%%%%%%%%%%%%%%%%%%%%%%%%%%%%%%%%%%%%%%%
\item{\bf c.\ } denoting by $d(\|Dp(u-\ell)\|)$ the 
measure in $E$ generated by the total variation 
$\|Dp(u-\ell)\|(E)$ 
%%%%%%%%%%%%%%%%%%%%%%%%%%%%%%%%%%%%%%%%%%%%%%%%%%%
\begin{equation}\label{Eq:4:2}
%%%%%%%%%%%%%%%%%%%%%%%%%%%%%%%%%%%%%%%%%%%%%%%%%%%
\begin{aligned}
%%%%%%%%%%%%%%%%%%%%%%%%%%%%%%%%%%%%%%%%%%%%%%%%%%%
&\ine\Big(\int_0^{u-\ell} p(s)ds\Big){\z(x,t_2)}\,dx 
+\int_{t_1}^{t_2}\ine\z d(\|D(p(u-\ell)\|)dt\\
&\le \ine \Big(\int_0^{u-\ell} p(s)ds\Big){\z(x,t_1)}\,dx 
-\int_{t_1}^{t_2}\ine \Big(\int_0^{u-\ell} p(s)ds\Big)\z_t dxdt\\
&\kern0.5cm -\int_{t_1}^{t_2}\ine{\bf z}\cdot D\z p(u-\ell)dxdt
%%%%%%%%%%%%%%%%%%%%%%%%%%%%%%%%%%%%%%%%%%%%%%%%%%%
\end{aligned}
%%%%%%%%%%%%%%%%%%%%%%%%%%%%%%%%%%%%%%%%%%%%%%%%%%%
\end{equation}
%%%%%%%%%%%%%%%%%%%%%%%%%%%%%%%%%%%%%%%%%%%%%%%%%%%
%%%%%%%%%%%%%%%%%%%%%%%%%%%%%%%%%%%%%%%%%%%%%%%%%%%
\end{description} %%%
%%%%%%%%%%%%%%%%%%%%%%%%%%%%%%%%%%%%%%%%%%%%%%%%%%%
for all $\ell\in\rr$, all $p\in\mcl{P}$, all 
$\z\in\mcl{C}(E_T)$ and all $[t_1,t_2]\subset(0,T)$. 
The notion is a local version of a global one 
introduced in \cite[Chapter~3]{ACM}. Similar notions are in 
\cite{ACM,BDM,BDS,LT}, associated with 
issues of existence for the Cauchy problem and 
boundary value problems associated with (\ref{Eq:4:1}).  
The notion of solution in \cite{BDM}, called {\it variational}, 
is different and closely related to the variational integrals 
(\ref{Eq:1:1}).  

Our results are local in nature and disengaged from any 
initial or boundary conditions. Let $u$ be a local solution 
to (\ref{Eq:4:1}) in the indicated sense, which in addition is 
locally bounded in $E_T$. In (\ref{Eq:4:2}) take $\ell=0$, 
and  $p_{\pm}(u)=\pm(u-k)_{\pm}$. Since $u\in L^\infty_{\loc}(E_T)$ 
one verifies that $p_{\pm}\in\mcl{P}$. Standard calculations 
then yield that $u$ is in the DeGiorgi classes $[DG]^{\pm}(E;\gm)$, 
for some fixed $\gm>0$.  As a consequence, we have the following:
%%%%%%%%%%%%%%%%%%%%%%%%%%%%%%%%%%%%%%%%%%%%%%%%%%%
\begin{corollary}\label{Cor:4:1}
%%%%%%%%%%%%%%%%%%%%%%%%%%%%%%%%%%%%%%%%%%%%%%%%%%%
Let $u\in L^\infty_{\loc}(E_T)$ be a local solution 
to (\ref{Eq:4:1}), in $E_T$, in the sense {\bf(a)-(c)} above. 
Then, $u$ is continuous at some $\pto\in E_T$, 
{if and only if} (\ref{Eq:1:3}) holds true.
%%%%%%%%%%%%%%%%%%%%%%%%%%%%%%%%%%%%%%%%%%%%%%%%%%%
\end{corollary}
%%%%%%%%%%%%%%%%%%%%%%%%%%%%%%%%%%%%%%%%%%%%%%%%%%%
{\begin{ack}
We  thank the referee for the valuable comments.
\end{ack}}
%%%%%%%%%%%%%%%%%%%%%%%%%%%%%%%%%%%%%%%%%%%%%%%%%%%%
\section{Proof of the Necessary Condition}\label{S:5}
%%%%%%%%%%%%%%%%%%%%%%%%%%%%%%%%%%%%%%%%%%%%%%%%%%%
Let $u\in [DG](E_T;\gm)$ be continuous at $\pto\in E_T$, 
which we may take as the origin of $\rr^{N+1}$, and may 
assume $u(0,0)=0$. In (\ref{Eq:3:2})$_+$ for $\ukp$, 
take $\theta=1$ and $k=0$. Let also $\z\in\mcl{C}(Q_{2\rho})$ be 
such that $\z(\cdot,-2\rho)=0$, such that 
$\z=1$ on $Q_{\frac32\rho}$, and 
%%%%%%%%%%%%%%%%%%%%%%%%%%%%%%%%%%%%%%%%%%%%%%%%%%%
\begin{equation*}
%%%%%%%%%%%%%%%%%%%%%%%%%%%%%%%%%%%%%%%%%%%%%%%%%%%
|D\z|+\z_t\le\frac3\rho.
%%%%%%%%%%%%%%%%%%%%%%%%%%%%%%%%%%%%%%%%%%%%%%%%%%%
\end{equation*}
%%%%%%%%%%%%%%%%%%%%%%%%%%%%%%%%%%%%%%%%%%%%%%%%%%%
Repeat the same choices in (\ref{Eq:3:2})${}_-$ 
for $\ukm$. Adding the resulting inequalities gives
%%%%%%%%%%%%%%%%%%%%%%%%%%%%%%%%%%%%%%%%%%%%%%%%%%%
\begin{equation}\label{Eq:5:1}
%%%%%%%%%%%%%%%%%%%%%%%%%%%%%%%%%%%%%%%%%%%%%%%%%%%
\frac\rho{|Q_\rho|}\int_{-2\rho}^0\|D (u\z)(\cdot,t)
\|(B_{2\rho})dt\le 2^{N+1}\gm 
\dashiint_{Q_{2\rho}}\big(u+u^2\big)dxdt.
%%%%%%%%%%%%%%%%%%%%%%%%%%%%%%%%%%%%%%%%%%%%%%%%%%%
\end{equation}
%%%%%%%%%%%%%%%%%%%%%%%%%%%%%%%%%%%%%%%%%%%%%%%%%%%
Since the total variation $\|Dw\|$ of a function $w\in BV$ 
can be seen as a measure (see, for example,  \cite[Chapter~1, 
\S~1]{ziemer}), we have
%%%%%%%%%%%%%%%%%%%%%%%%%%%%%%%%%%%%%%%%%%%%%%%%%%%
\begin{equation*}
%%%%%%%%%%%%%%%%%%%%%%%%%%%%%%%%%%%%%%%%%%%%%%%%%%%
\frac\rho{|Q_\rho|}\int_{-\rho}^0\|D (u\z)(\cdot,t)
\|(B_{\rho})dt\le\frac\rho{|Q_\rho|}\int_{-2\rho}^0
\|D (u\z)(\cdot,t)\|(B_{2\rho})dt;
%%%%%%%%%%%%%%%%%%%%%%%%%%%%%%%%%%%%%%%%%%%%%%%%%%%
\end{equation*}
%%%%%%%%%%%%%%%%%%%%%%%%%%%%%%%%%%%%%%%%%%%%%%%%%%%
on the other hand, $u\z\equiv u$ in $Q_{\frac32\rho}\supset Q_\rho$, 
and therefore we conclude
%%%%%%%%%%%%%%%%%%%%%%%%%%%%%%%%%%%%%%%%%%%%%%%%%%%
\begin{equation*}
%%%%%%%%%%%%%%%%%%%%%%%%%%%%%%%%%%%%%%%%%%%%%%%%%%%
\frac\rho{|Q_\rho|}\int_{-\rho}^0\|D u(\cdot,t)\|(B_{\rho})dt\le 2^{N+1}\gm\, 
\dashiint_{Q_{2\rho}}\big(u+u^2\big)dxdt.
%%%%%%%%%%%%%%%%%%%%%%%%%%%%%%%%%%%%%%%%%%%%%%%%%%%
\end{equation*}
%%%%%%%%%%%%%%%%%%%%%%%%%%%%%%%%%%%%%%%%%%%%%%%%%%%
The right-hand side tends to zero as $\rho\to0$, thereby 
implying the necessary condition of Theorem~\ref{Thm:1:1}.\hfill\bbox
%%%%%%%%%%%%%%%%%%%%%%%%%%%%%%%%%%%%%%%%%%%%%%%%%
%%%%%%%%%%%%%%%%%%%%%%%%%%%%%%%%%%%%%%%%%%%%%%
\section{A DeGiorgi-Type Lemma}\label{S:6}
%%%%%%%%%%%%%%%%%%%%%%%%%%%%%%%%%%%%%%%%%%%%%%
For a fixed cylinder $[(y,s)+Q_{2\rho}(\theta)]\subset E_T$, 
 denote by $\mu_{\pm}$ and $\om$, non-negative numbers 
such that
%%%%%%%%%%%%%%%%%%%%%%%%%%%%%%%%%%%%%%%%%%%%%%
\begin{equation}\label{Eq:6:1}
%%%%%%%%%%%%%%%%%%%%%%%%%%%%%%%%%%%%%%%%%%%%%%
\mu_+\ge \essup_{[(y,s)+Q_{2\rho}(\theta)]} u,\quad 
\mu_-\le \essinf_{[(y,s)+Q_{2\rho}(\theta)]} u,\quad 
\om\ge\mu_+-\mu_-.
%%%%%%%%%%%%%%%%%%%%%%%%%%%%%%%%%%%%%%%%%%%%%%
\end{equation}
%%%%%%%%%%%%%%%%%%%%%%%%%%%%%%%%%%%%%%%%%%%%%%
Let $\xi\in(0,\frac12]$ be fixed and let $\theta=2\xi\om$. 
{This is an intrinsic cylinder in that its length 
$\theta\rho$ depends on the oscillation of $u$ within it. 
We assume momentarily that the indicated choice of 
parameters can be effected.}
%%%%%%%%%%%%%%%%%%%%%%%%%%%%%%%%%%%%%%%%%%%%%%%%%%%%%%%%%%%%%%%%
\begin{lemma}\label{Lm:6:1} 
%%%%%%%%%%%%%%%%%%%%%%%%%%%%%%%%%%%%%%%%%%%%%%
Let $u$ belong to  $[DG]^-(E_T,\gm)$. 
There exists a number $\nu_-$ depending on 
$N$, and $\gm$ only, such that if 
%%%%%%%%%%%%%%%%%%%%%%%%%%%%%%%%%%%%%%%%%%%%%%
\begin{equation}\label{Eq:6:2}
%%%%%%%%%%%%%%%%%%%%%%%%%%%%%%%%%%%%%%%%%%%%%%
\left|[u\le\mu_-+\xi\om]\cap 
[(y,s)+\qrtt]\right|\le 
\nu_-|\qrtt|,
%%%%%%%%%%%%%%%%%%%%%%%%%%%%%%%%%%%%%%%%%%%%%%
\end{equation}
%%%%%%%%%%%%%%%%%%%%%%%%%%%%%%%%%%%%%%%%%%%%%%
then 
%%%%%%%%%%%%%%%%%%%%%%%%%%%%%%%%%%%%%%%%%%%%%%
\begin{equation}\label{Eq:6:3}
%%%%%%%%%%%%%%%%%%%%%%%%%%%%%%%%%%%%%%%%%%%%%%
u\ge \mu_-+{\txty\frac12}\xi\om\quad\text{ a.e. in  }\>\big[
(y,s)+\qrt\big].
%%%%%%%%%%%%%%%%%%%%%%%%%%%%%%%%%%%%%%%%%%%%%%
\end{equation}
%%%%%%%%%%%%%%%%%%%%%%%%%%%%%%%%%%%%%%%%%%%%%%
Likewise, if $u$ belongs to  $[DG]^+(E_T,\gm)$, 
there exists a number $\nu_+$ depending on 
$N$, and $\gm$ only, such that if
%%%%%%%%%%%%%%%%%%%%%%%%%%%%%%%%%%%%%%%%%%%%%%
\begin{equation}\label{Eq:6:4}
%%%%%%%%%%%%%%%%%%%%%%%%%%%%%%%%%%%%%%%%%%%%%%
\left|[u\ge \mu_+-\xi\om]\cap 
[(y,s)+\qrtt]\right|\le 
\nu_+|\qrtt|,
%%%%%%%%%%%%%%%%%%%%%%%%%%%%%%%%%%%%%%%%%%%%%%
\end{equation}
%%%%%%%%%%%%%%%%%%%%%%%%%%%%%%%%%%%%%%%%%%%%%%
then 
%%%%%%%%%%%%%%%%%%%%%%%%%%%%%%%%%%%%%%%%%%%%%%
\begin{equation}\label{Eq:6:5}
%%%%%%%%%%%%%%%%%%%%%%%%%%%%%%%%%%%%%%%%%%%%%%
u\le\mu_+-{\txty\frac12}\xi\om\quad\text{ a.e. in  }\>\big[
(y,s)+\qrt\big].
%%%%%%%%%%%%%%%%%%%%%%%%%%%%%%%%%%%%%%%%%%%%%%
\end{equation}
%%%%%%%%%%%%%%%%%%%%%%%%%%%%%%%%%%%%%%%%%%%%%%
\end{lemma}
%%%%%%%%%%%%%%%%%%%%%%%%%%%%%%%%%%%%%%%%%%%%%%
%%%%%%%%%%%%%%%%%%%%%%%%%%%%%%%%%%%%%%%%%%%%%%%%%%%
%%%%%%%%%%%%%%%%%%%%%%%%%%%%%%%%%%%%%%%%%%%%%%%%%%%
\noi{\bf Proof: } We prove (\ref{Eq:6:2})--(\ref{Eq:6:3}), the 
proof for \eqref{Eq:6:4}--\eqref{Eq:6:5} being similar. 
We may assume $(y,s)=(0,0)$ and for $n=0,1,\dots$, set
%%%%%%%%%%%%%%%%%%%%%%%%%%%%%%%%%%%%%%%%%%%%%%%%%%%
\begin{equation*}%%\label{Eq:6:5}
%%%%%%%%%%%%%%%%%%%%%%%%%%%%%%%%%%%%%%%%%%%%%%%%%%%
\rho_n=\rho+\frac{\rho}{2^{n}},
\qquad B_n=B_{\rho_n},\qquad Q_n=B_n\times(-\theta\rho_n,0].
%%%%%%%%%%%%%%%%%%%%%%%%%%%%%%%%%%%%%%%%%%%%%%%%%%%
\end{equation*}
%%%%%%%%%%%%%%%%%%%%%%%%%%%%%%%%%%%%%%%%%%%%%%%%%%%
Apply \eqref{Eq:3:2}$_-$ over $B_n$ and $Q_n$ to $\uknm$, 
for the levels 
%%%%%%%%%%%%%%%%%%%%%%%%%%%%%%%%%%%%%%%%%%%%%%%%%%%
\begin{equation*}%%\label{Eq:6:6}
%%%%%%%%%%%%%%%%%%%%%%%%%%%%%%%%%%%%%%%%%%%%%%%%%%%
k_n=\mu_-+\xi_n\om\qquad\text{ where }\qquad 
\xi_n=\frac12\xi+\frac1{2^{n+1}}\xi.
%%%%%%%%%%%%%%%%%%%%%%%%%%%%%%%%%%%%%%%%%%%%%%%%%%%
\end{equation*}
%%%%%%%%%%%%%%%%%%%%%%%%%%%%%%%%%%%%%%%%%%%%%%%%%%%
The cutoff function $\z$ is taken of the form 
$\z(x,t)=\z_1(x)\z_2(t)$, where 
%%%%%%%%%%%%%%%%%%%%%%%%%%%%%%%%%%%%%%%%%%%%%%%%%%%
\begin{equation*}%%\label{Eq:6:7}
%%%%%%%%%%%%%%%%%%%%%%%%%%%%%%%%%%%%%%%%%%%%%%%%%%%
\begin{array}{lc}
%%%%%%%%%%%%%%%%%%%%%%%%%%%%%%%%%%%%%%%%%%%%%%%%%%%
{\dsty 
\z_1=\left\{
\begin{array}{ll}
1\>&\text{ in }\> B_{n+1}\\
{}\\
0\>&\text{ in }\>\rn-B_n
\end{array}\right .}\quad
&{\dsty |D\z_1|\le \frac1{\rho_n-\rho_{n+1}}
=\frac{2^{n+1}}{\rho} }\\
{}\\
{\dsty 
\z_2=\left\{
\begin{array}{ll}
0\>&\text{ for }\> t<-\theta \rho_n\\
{}\\
1\>&\text{ for }\> t\ge-\theta\rho_{n+1}
\end{array}\right .}\quad
&{\dsty 0\le \z_{2,t}\le\frac1{\theta(\rho_n-\rho_{n+1})}
=\frac{2^{(n+1)}}{\theta \rho}}.
%%%%%%%%%%%%%%%%%%%%%%%%%%%%%%%%%%%%%%%%%%%%%%%%%%%
\end{array}
%%%%%%%%%%%%%%%%%%%%%%%%%%%%%%%%%%%%%%%%%%%%%%%%%%%
\end{equation*} 
%%%%%%%%%%%%%%%%%%%%%%%%%%%%%%%%%%%%%%%%%%%%%%%%%%%
Inequality (\ref{Eq:3:2})$_-$ with these 
stipulations yields
%%%%%%%%%%%%%%%%%%%%%%%%%%%%%%%%%%%%%%%%%%%%%%%%%%%
\begin{equation*}%%\label{Eq:6:8}
%%%%%%%%%%%%%%%%%%%%%%%%%%%%%%%%%%%%%%%%%%%%%%%%%%%
\begin{aligned}
%%%%%%%%%%%%%%%%%%%%%%%%%%%%%%%%%%%%%%%%%%%%%%%%%%%
&\essup_{-\theta\rho_n<t<0}\int_{B_n}\uknm^2\z(x,t)dx 
+\int_{-\theta\rho_n}^0\|D\uknm\z\|(B_n)dt\\
&\le\gm\frac{2^{n}}{\rho}\left(
\iint_{Q_n}\uknm dxdt+\frac{1}{\theta}
\iint_{Q_n}\uknm^2 dxdt\right)\\
&\le\gm\frac{2^{n}(\xi\om)}{\rho}|[u<k_n]\cap Q_n|.
%%%%%%%%%%%%%%%%%%%%%%%%%%%%%%%%%%%%%%%%%%%%%%%%%%%
\end{aligned}
%%%%%%%%%%%%%%%%%%%%%%%%%%%%%%%%%%%%%%%%%%%%%%%%%%%
\end{equation*}
%%%%%%%%%%%%%%%%%%%%%%%%%%%%%%%%%%%%%%%%%%%%%%%%%%%
By the embedding Proposition~{4.1} of 
\cite[Preliminaries]{DBGV-mono}
%%%%%%%%%%%%%%%%%%%%%%%%%%%%%%%%%%%%%%%%%%%%%%%%%%%
\begin{equation*}
%%%%%%%%%%%%%%%%%%%%%%%%%%%%%%%%%%%%%%%%%%%%%%%%%%%
\begin{aligned}
%%%%%%%%%%%%%%%%%%%%%%%%%%%%%%%%%%%%%%%%%%%%%%%%%%%
\iint_{Q_n}&[\uknm\z]^{\frac{N+2}N} dxdt
\le\int_{-\theta\rho_n}^0\|D[\uknm\z]\|(B_n)dt\\
&\qquad\times\left(\essup_{-\theta\rho_n<t<0} 
\int_{B_n}[\uknm\z(x,t)]^2dx\right)^{\frac{1}{N}}\\
&\le\gm\left(\frac{2^{n}}{\rho}\xi\om
\right)^{\frac{N+1}N}|[u<k_n]\cap Q_n|^{\frac{N+1}N}.
%%%%%%%%%%%%%%%%%%%%%%%%%%%%%%%%%%%%%%%%%%%%%%%%%%%
\end{aligned}
%%%%%%%%%%%%%%%%%%%%%%%%%%%%%%%%%%%%%%%%%%%%%%%%%%%
\end{equation*}
%%%%%%%%%%%%%%%%%%%%%%%%%%%%%%%%%%%%%%%%%%%%%%%%%%%
Estimate below 
%%%%%%%%%%%%%%%%%%%%%%%%%%%%%%%%%%%%%%%%%%%%%%%%%%%
\begin{equation*}
%%%%%%%%%%%%%%%%%%%%%%%%%%%%%%%%%%%%%%%%%%%%%%%%%%%
\iint_{Q_n}[\uknm\z]^{\frac{N+2}N} dxdt
\ge\left(\frac{\xi\om}{2^{n+2}}\right)^{\frac{N+2}N}
|[u<k_{n+1}]\cap Q_{n+1}|
%%%%%%%%%%%%%%%%%%%%%%%%%%%%%%%%%%%%%%%%%%%%%%%%%%%
\end{equation*}
%%%%%%%%%%%%%%%%%%%%%%%%%%%%%%%%%%%%%%%%%%%%%%%%%%%
and set 
%%%%%%%%%%%%%%%%%%%%%%%%%%%%%%%%%%%%%%%%%%%%%%%%%%%
\begin{equation*}
%%%%%%%%%%%%%%%%%%%%%%%%%%%%%%%%%%%%%%%%%%%%%%%%%%%
Y_n=\frac{|[u<k_n]\cap Q_n|}{|Q_n|}.
%%%%%%%%%%%%%%%%%%%%%%%%%%%%%%%%%%%%%%%%%%%%%%%%%%%
\end{equation*}
%%%%%%%%%%%%%%%%%%%%%%%%%%%%%%%%%%%%%%%%%%%%%%%%%%%
Then
%%%%%%%%%%%%%%%%%%%%%%%%%%%%%%%%%%%%%%%%%%%%%%%%%%%
\begin{equation*}%%
%%%%%%%%%%%%%%%%%%%%%%%%%%%%%%%%%%%%%%%%%%%%%%%%%%%
Y_{n+1}\le\gm b^nY_n^{1+\frac{1}{N}}
%%%%%%%%%%%%%%%%%%%%%%%%%%%%%%%%%%%%%%%%%%%%%%%%%%%
\end{equation*} 
%%%%%%%%%%%%%%%%%%%%%%%%%%%%%%%%%%%%%%%%%%%%%%%%%%%
where
%%%%%%%%%%%%%%%%%%%%%%%%%%%%%%%%%%%%%%%%%%%%%%%%%%%
\begin{equation*}%%\label{Eq:6:9}
%%%%%%%%%%%%%%%%%%%%%%%%%%%%%%%%%%%%%%%%%%%%%%%%%%%
b=2^{\frac 1N[3N+4]}.
%%%%%%%%%%%%%%%%%%%%%%%%%%%%%%%%%%%%%%%%%%%%%%%%%%%
\end{equation*}
%%%%%%%%%%%%%%%%%%%%%%%%%%%%%%%%%%%%%%%%%%%%%%%%%%%
By Lemma~{5.1} of \cite[Preliminaries]{DBGV-mono}, $\{Y_n\}\to0$ 
as $n\to\infty$, provided 
%%%%%%%%%%%%%%%%%%%%%%%%%%%%%%%%%%%%%%%%%%%%%%%%%%%
\begin{equation*}%%\label{Eq:6:10}
%%%%%%%%%%%%%%%%%%%%%%%%%%%%%%%%%%%%%%%%%%%%%%%%%%%
Y_o\le\gm^{-N}b^{-{N}^2}\,\df{=}\, \nu_-.
%%%%%%%%%%%%%%%%%%%%%%%%%%%%%%%%%%%%%%%%%%%%%%%%%%%
\end{equation*}
%%%%%%%%%%%%%%%%%%%%%%%%%%%%%%%%%%%%%%%%%%%%%%%%%%%%%%%%%%%%%%%%%%%%
The proof of (\ref{Eq:6:4})--(\ref{Eq:6:5}) is 
almost identical. One starts from inequalities 
(\ref{Eq:3:2})$_+$ written for the truncated functions
%%%%%%%%%%%%%%%%%%%%%%%%%%%%%%%%%%%%%%%%%%%%%%%%%%%%%%%%%%%%%%%%%%%%
\begin{equation*}
%%%%%%%%%%%%%%%%%%%%%%%%%%%%%%%%%%%%%%%%%%%%%%%%%%%%%%%%%%%%%%%%%%%%
\uknp\qquad\text{ with }\qquad k_n=\mu_+-\xi_n\om
%%%%%%%%%%%%%%%%%%%%%%%%%%%%%%%%%%%%%%%%%%%%%%%%%%%%%%%%%%%%%%%%%%%%
\end{equation*} 
%%%%%%%%%%%%%%%%%%%%%%%%%%%%%%%%%%%%%%%%%%%%%%%%%%%%%%%%%%%%%%%%%%%%
for the same choice of $\xi_n$.\hfill\bbox
%%%%%%%%%%%%%%%%%%%%%%%%%%%%%%%%%%%%%%%%%%%%%%%%%%%%%
%%%%%%%%%%%%%%%%%%%%%%%%%%%%%%%%%%%%%
\section{A Time Expansion of Positivity}\label{S:7}
%%%%%%%%%%%%%%%%%%%%%%%%%%%%%%%%%%%%%%%%%%%%%%%%%%%%%%%
For a fixed cylinder  
%%%%%%%%%%%%%%%%%%%%%%%%%%%%%%%%%%%%%%%%%%%%%%%%%%%%%%%
\begin{equation*}
%%%%%%%%%%%%%%%%%%%%%%%%%%%%%%%%%%%%%%%%%%%%%%%%%%%%%%%
[(y,s)+Q_{2\rho}^+(\theta)]=B_{2\rho}(y)
\times(s,s+\theta\rho) \subset E_T, 
%%%%%%%%%%%%%%%%%%%%%%%%%%%%%%%%%%%%%%%%%%%%%%%%%%%%%%%
\end{equation*}
%%%%%%%%%%%%%%%%%%%%%%%%%%%%%%%%%%%%%%%%%%%%%%%%%%%%%%%
denote by $\mu_{\pm}$ and $\om$, non-negative numbers 
satisfying the analog of (\ref{Eq:6:1}).
%%%%%%%%%%%%%%%%%%%%%%%%%%%%%%%%%%%%%%%%%%%%%%%%%%%%%%%
Let also $\xi\in(0,1)$ be a fixed parameter. {The value of $\theta$ will be determined by the proof; we momentarily assume that such a choice can be done.}
%%%%%%%%%%%%%%%%%%%%%%%%%%%%%%%%%
\begin{lemma}\label{Lm:7:1} 
%%%%%%%%%%%%%%%%%%%%%%%%%%%%%%%%%%%%%%%%%%%%%%%%%%%%%%%
Let $u\in [DG]^-(E_T,\gm)$ and assume that 
for some $(y,s)\in E_T$ and some $\rho>0$  
%%%%%%%%%%%%%%%%%%%%%%%%%%%%%%%%%%%%%%%%%%%%%%%%%%%%%%%
\begin{equation*}
%%%%%%%%%%%%%%%%%%%%%%%%%%%%%%%%%%%%%%%%%%%%%%%%%%%%%%%
\big|[u(\cdot,s)\ge \mu_-+\xi\om]\cap B_\rho(y)\big| 
\ge{\txty\frac12}\big|B_{\rho}(y)\big|. 
%%%%%%%%%%%%%%%%%%%%%%%%%%%%%%%%%%%%%%%%%%%%%%%%%%%%%%%
\end{equation*}
%%%%%%%%%%%%%%%%%%%%%%%%%%%%%%%%%%%%%%%%%%%%%%%%%%%%%%%
Then, there exist $\dl$ and $\eps$ in $(0,1)$, depending 
only on $N$, $\gm$, and independent of $\xi$, such that 
%%%%%%%%%%%%%%%%%%%%%%%%%%%%%%%%%%%%%%%%%%%%%%%%%%%%%%%
\begin{equation*}%%\label{Eq:7:1}
%%%%%%%%%%%%%%%%%%%%%%%%%%%%%%%%%%%%%%%%%%%%%%%%%%%%%%%
\big|[u(\cdot,t)>\mu_-+\eps\xi\om]\cap B_\rho(y)\big|
\ge{\txty\frac14}|B_\rho|\quad\text{ for all }\> 
t \in\big(s,s+{\dl (\xi\om)\rho}\big].
%%%%%%%%%%%%%%%%%%%%%%%%%%%%%%%%%%%%%%%%%%%%%%%%%%%%%%%
\end{equation*}
%%%%%%%%%%%%%%%%%%%%%%%%%%%%%%%%%%%%%%%%%%%%%%%%%%%%%%%
\end{lemma}
%%%%%%%%%%%%%%%%%%%%%%%%%%%%%%%%%%%%%%%%%%%%%%%%%%%%%%%%%%%%%%%%%%%%%
\noi{\bf Proof: } Assume $(y,s)=(0,0)$ and for 
$k>0$ and $t>0$ set 
%%%%%%%%%%%%%%%%%%%%%%%%%%%%%%%%%%%%%%%%%%%%%%%%%%%%%%%%%%%%%%%%%%%%%
\begin{equation*}%%\label{Eq:2:1:3}
%%%%%%%%%%%%%%%%%%%%%%%%%%%%%%%%%%%%%%%%%%%%%%%%%%%%%%%%%%%%%%%%%%%%%
A_{k,\rho}(t)=[u(\cdot,t)<k]\cap B_{\rho}.
%%%%%%%%%%%%%%%%%%%%%%%%%%%%%%%%%%%%%%%%%%%%%%%%%%%%%%%%%%%%%%%%%%%%%
\end{equation*}
%%%%%%%%%%%%%%%%%%%%%%%%%%%%%%%%%%%%%%%%%%%%%%%%%%%%%%%%%%%%%%%%%%%%%
The assumption implies
%%%%%%%%%%%%%%%%%%%%%%%%%%%%%%%%%%%%%%%%%%%%%%%%%%%%%%%%%%%%%%%%%%%%%
\begin{equation}\label{Eq:7:1}
%%%%%%%%%%%%%%%%%%%%%%%%%%%%%%%%%%%%%%%%%%%%%%%%%%%%%%%%%%%%%%%%%%%%%
|A_{\mu_-+\xi\om,\rho}(0)|\le{\txty\frac12}|B_\rho|.
%%%%%%%%%%%%%%%%%%%%%%%%%%%%%%%%%%%%%%%%%%%%%%%%%%%%%%%%%%%%%%%%%%%%%
\end{equation}
%%%%%%%%%%%%%%%%%%%%%%%%%%%%%%%%%%%%%%%%%%%%%%%%%%%%%%%%%%%%%%%%%%%%%
Write down inequalities \eqref{Eq:3:2}$_-$
for the truncated functions $(u-(\mu_-+\xi\om))_-$, over 
the cylinder $B_\rho\times(0,\theta\rho]$, where 
$\theta>0$ is to be chosen. The cutoff function $\z$ 
is taken independent of $t$, non-negative, and such that
%%%%%%%%%%%%%%%%%%%%%%%%%%%%%%%%%%%%%%%%%%%%%%%%%%%%%%%%%%%%%%%%%%%%%
\begin{equation*}
%%%%%%%%%%%%%%%%%%%%%%%%%%%%%%%%%%%%%%%%%%%%%%%%%%%%%%%%%%%%%%%%%%%%%
\z=1 \quad\text{ on }\> B_{(1-\sig)\rho},
\quad\text{ and }\quad 
|D\z|\le\frac1{\sig\rho},
%%%%%%%%%%%%%%%%%%%%%%%%%%%%%%%%%%%%%%%%%%%%%%%%%%%%%%%%%%%%%%%%%%%%%
\end{equation*}
%%%%%%%%%%%%%%%%%%%%%%%%%%%%%%%%%%%%%%%%%%%%%%%%%%%%%%%%%%%%%%%%%%%%%
where $\sig\in(0,1)$ is to be chosen. Discarding  the 
non-negative term containing $D(u-(\mu_-+\xi\om))_-$ 
on the left-hand side, these inequalities yield
%%%%%%%%%%%%%%%%%%%%%%%%%%%%%%%%%%%%%%%%%%%%%%%%%%%%%%%%%%%%%%%%%%%%%
\begin{align*}
%%%%%%%%%%%%%%%%%%%%%%%%%%%%%%%%%%%%%%%%%%%%%%%%%%%%%%%%%%%%%%%%%%%%%
\int_{B_{(1-\sig)\rho}}&(u-(\mu_-+\xi\om))_-^2(x,t)dx\le
\int_{B_\rho}(u-(\mu_-+\xi\om))_-^2(x,0)dx\\
&\quad+\frac{\gm}{\sig\rho}\int_0^{\theta\rho}
\int_{B_\rho}(u-(\mu_-+\xi\om))_- dxdt\\
&\le (\xi\om)^2\Big[\frac12+\gm\frac{\theta}{\sig(\xi\om)}\Big]
|B_\rho|
%%%%%%%%%%%%%%%%%%%%%%%%%%%%%%%%%%%%%%%%%%%%%%%%%%%%%%%%%%%%%%%%%%%%%
\end{align*}
%%%%%%%%%%%%%%%%%%%%%%%%%%%%%%%%%%%%%%%%%%%%%%%%%%%%%%%%%%%%%%%%%%%%%
for all $t\in(0,\theta\rho]$, where we have 
enforced (\ref{Eq:7:1}).  The left-hand side is 
estimated below by
%%%%%%%%%%%%%%%%%%%%%%%%%%%%%%%%%%%%%%%%%%%%%%%%%%%%%%%%%%%%%%%%%%%%%
\begin{align*}
%%%%%%%%%%%%%%%%%%%%%%%%%%%%%%%%%%%%%%%%%%%%%%%%%%%%%%%%%%%%%%%%%%%%%
&\int_{B_{(1-\sig)\rho}}(u-(\mu_-+\xi\om))_-^2(x,t)dx\\
&\ge
\int_{B_{(1-\sig)\rho}\cap [u<\mu_-+\eps\xi\om]}(u-(\mu_-+\xi\om))_-^2(x,t)dx\\
&\ge (\xi\om)^2(1-\eps)^2|A_{\mu_-+\eps\xi\om,(1-\sig)\rho}(t)|
%%%%%%%%%%%%%%%%%%%%%%%%%%%%%%%%%%%%%%%%%%%%%%%%%%%%%%%%%%%%%%%%%%%%%
\end{align*}
%%%%%%%%%%%%%%%%%%%%%%%%%%%%%%%%%%%%%%%%%%%%%%%%%%%%%%%%%%%%%%%%%%%%%
where $\eps\in(0,1)$ is to be chosen. Next, estimate
%%%%%%%%%%%%%%%%%%%%%%%%%%%%%%%%%%%%%%%%%%%%%%%%%%%%%%%%%%%%%%%%%%%%%
\begin{align*}
%%%%%%%%%%%%%%%%%%%%%%%%%%%%%%%%%%%%%%%%%%%%%%%%%%%%%%%%%%%%%%%%%%%%%
|A_{\mu_-+\eps\xi\om,\rho}(t)|&=|A_{\mu_-+\eps\xi\om,(1-\sig)\rho}(t)\cup 
(A_{\mu_-+\eps\xi\om,\rho}(t)-A_{\mu_-+\eps\xi\om,(1-\sig)\rho}(t))|\\
&\le|A_{\mu_-+\eps\xi\om,(1-\sig)\rho}(t)|+|B_\rho-B_{(1-\sig)\rho}|\\
&\le |A_{\mu_-+\eps\xi\om,(1-\sig)\rho}(t)|+N\sig|B_\rho|.
%%%%%%%%%%%%%%%%%%%%%%%%%%%%%%%%%%%%%%%%%%%%%%%%%%%%%%%%%%%%%%%%%%%%%
\end{align*}
%%%%%%%%%%%%%%%%%%%%%%%%%%%%%%%%%%%%%%%%%%%%%%%%%%%%%%%%%%%%%%%%%%%%%
Combining these estimates gives
%%%%%%%%%%%%%%%%%%%%%%%%%%%%%%%%%%%%%%%%%%%%%%%%%%%%%%%%%%%%%%%%%%%%%
\begin{align*}
%%%%%%%%%%%%%%%%%%%%%%%%%%%%%%%%%%%%%%%%%%%%%%%%%%%%%%%%%%%%%%%%%%%%%
|A_{\mu_-+\eps\xi\om,\rho}(t)|&\le\frac1{(\xi\om)^2(1-\eps)^2}
\int_{B_{(1-\sig)\rho}}(u-(\mu_-+\xi\om))_-^2(x,t)dx+N\sig|B_{\rho}|\\
&\le\frac1{(1-\eps)^2}\Big[\frac12+
\frac{\gm\theta}{\sig(\xi\om)}+N\sig\Big]|B_{\rho}|.
%%%%%%%%%%%%%%%%%%%%%%%%%%%%%%%%%%%%%%%%%%%%%%%%%%%%%%%%%%%%%%%%%%%%%
\end{align*}
%%%%%%%%%%%%%%%%%%%%%%%%%%%%%%%%%%%%%%%%%%%%%%%%%%%%%%%%%%%%%%%%%%%%%
Choose $\theta=\dl(\xi\om)$ and then set 
%%%%%%%%%%%%%%%%%%%%%%%%%%%%%%%%%%%%%%%%%%%%%%%%%%%%%%%%%%%%%%%%%%%%%
\begin{equation}\label{Eq:7:2}
%%%%%%%%%%%%%%%%%%%%%%%%%%%%%%%%%%%%%%%%%%%%%%%%%%%%%%%%%%%%%%%%%%%%%
\sig=\frac{1}{16N},\qquad 
\eps\le\frac1{32},\qquad \dl=\frac{1}{2^8\gm N}.
%%%%%%%%%%%%%%%%%%%%%%%%%%%%%%%%%%%%%%%%%%%%%%%%%%%%%%%%%%%%%%%%%%%%%
\end{equation} 
%%%%%%%%%%%%%%%%%%%%%%%%%%%%%%%%%%%%%%%%%%%%%%%%%%%%%%%%%%%%%%%%%%%%%
This proves the lemma.\hfill\bbox
%%%%%%%%%%%%%%%%%%%%%%%%%%%%%%%%%%%%%%%%%%%%%%%%%%
\section{Proof of the Sufficient Part of 
Theorem~\ref{Thm:1:1}}\label{S:8}
%%%%%%%%%%%%%%%%%%%%%%%%%%%%%%%%%%%%%%%%%%%%%%%%%%
Having fixed $\pto\in E_T$ assume it coincides 
with the origin of $\rr^{N+1}$ and let $\rho>0$ 
be so small that $Q_{\rho}\subset E_T$. Set 
%%%%%%%%%%%%%%%%%%%%%%%%%%%%%%%%%%%%%%%%%%%%%%%%%%
\begin{equation*}%%\label{Eq:6:1}
%%%%%%%%%%%%%%%%%%%%%%%%%%%%%%%%%%%%%%%%%%%%%%%%%%
\mu_+=\essup_{Q_{\rho}}u,\quad 
\mu_-=\essinf_{Q_{\rho}}u,\quad 
\om=\mu_+-\mu_-=\essosc_{Q_{\rho}} u.
%%%%%%%%%%%%%%%%%%%%%%%%%%%%%%%%%%%%%%%%%%%%%%%%%%
\end{equation*}
%%%%%%%%%%%%%%%%%%%%%%%%%%%%%%%%%%%%%%%%%%%%%%%%%%
Without loss of generality, we may assume 
that $\om\le1$ so that 
%%%%%%%%%%%%%%%%%%%%%%%%%%%%%%%%%%%%%%%%%%%%%%%%%%
\begin{equation*}%%\label{Eq:6:2}
%%%%%%%%%%%%%%%%%%%%%%%%%%%%%%%%%%%%%%%%%%%%%%%%%%
Q_\rho(\om)=B_{\rho}\times
(-\om\rho,0]\subset Q_{\rho}\subset E_T
%%%%%%%%%%%%%%%%%%%%%%%%%%%%%%%%%%%%%%%%%%%%%%%%%%
\end{equation*}
%%%%%%%%%%%%%%%%%%%%%%%%%%%%%%%%%%%%%%%%%%%%%%%%%%
and 
%%%%%%%%%%%%%%%%%%%%%%%%%%%%%%%%%%%%%%%%%%%%%%%%%%
\begin{equation*}
%%%%%%%%%%%%%%%%%%%%%%%%%%%%%%%%%%%%%%%%%%%%%%%%%%
\essosc_{Q_\rho(\om)} u\le \om.
%%%%%%%%%%%%%%%%%%%%%%%%%%%%%%%%%%%%%%%%%%%%%%%%%%
\end{equation*}
%%%%%%%%%%%%%%%%%%%%%%%%%%%%%%%%%%%%%%%%%%%%%%%%%%
If $u$ were not continuous at $\pto$, there would 
exist $\rho_o>0$ and $\om_o>0$, such that 
%%%%%%%%%%%%%%%%%%%%%%%%%%%%%%%%%%%%%%%%%%%%%%%%%%
\begin{equation}\label{Eq:8:1}
%%%%%%%%%%%%%%%%%%%%%%%%%%%%%%%%%%%%%%%%%%%%%%%%%%
\om_\rho=\essosc_{Q_{\rho}} u\ge \om_o>0
\quad\text{ for all }\> \rho\le\rho_o. 
%%%%%%%%%%%%%%%%%%%%%%%%%%%%%%%%%%%%%%%%%%%%%%%%%%
\end{equation}
%%%%%%%%%%%%%%%%%%%%%%%%%%%%%%%%%%%%%%%%%%%%%%%%%%
Let $\dl$ be determined from the last of (\ref{Eq:7:2}).  
 At the time level $t=-\dl\om\rho$, either 
%%%%%%%%%%%%%%%%%%%%%%%%%%%%%%%%%%%%%%%%%%%%%%%%%%
\begin{equation*}
%%%%%%%%%%%%%%%%%%%%%%%%%%%%%%%%%%%%%%%%%%%%%%%%%%
\begin{aligned}
%%%%%%%%%%%%%%%%%%%%%%%%%%%%%%%%%%%%%%%%%%%%%%%%%%
&\big|\big[u(\cdot,-\dl\om\rho)\ge 
\mu_-+{\txty\frac12}\om\big]\cap B_\rho\big|
\ge{\txty\frac12} |B_{\rho}|,\quad\text{ or }\\
%%%%%%%%%%%%%%%%%%%%%%%%%%%%%%%%%%%%%%%%%%%%%%%%%%
&\big|\big[u(\cdot,-\dl\om\rho)\le 
\mu_+-{\txty\frac12}\om\big]\cap B_\rho\big| 
\ge{\txty\frac12} |B_{\rho}|.
%%%%%%%%%%%%%%%%%%%%%%%%%%%%%%%%%%%%%%%%%%%%%%%%%%
\end{aligned}
%%%%%%%%%%%%%%%%%%%%%%%%%%%%%%%%%%%%%%%%%%%%%%%%%%
%%%%%%%%%%%%%%%%%%%%%%%%%%%%%%%%%%%%%%%%%%%%%%%%%%
\end{equation*}
%%%%%%%%%%%%%%%%%%%%%%%%%%%%%%%%%%%%%%%%%%%%%%%%%%
Assuming the former holds, by Lemma~\ref{Lm:7:1} 
%%%%%%%%%%%%%%%%%%%%%%%%%%%%%%%%%%%%%%%%%%%%%%%%%%
\begin{equation*}%%\label{Eq:6:2}
%%%%%%%%%%%%%%%%%%%%%%%%%%%%%%%%%%%%%%%%%%%%%%%%%%
\big|\big[u(\cdot,t)>\mu_-+{\txty\frac1{64}}\om\big]\cap B_\rho\big|
\ge{\txty\frac14}|B_\rho|\quad\text{ for all }\> 
t \in(-\dl\om\rho,0].
%%%%%%%%%%%%%%%%%%%%%%%%%%%%%%%%%%%%%%%%%%%%%%%%%%
\end{equation*}
%%%%%%%%%%%%%%%%%%%%%%%%%%%%%%%%%%%%%%%%%%%%%%%%%%
Let $2\xi=\frac1{64}\dl$. Then
%%%%%%%%%%%%%%%%%%%%%%%%%%%%%%%%%%%%%%%%%%%%%%%%%%
\begin{equation}\label{Eq:8:2}
%%%%%%%%%%%%%%%%%%%%%%%%%%%%%%%%%%%%%%%%%%%%%%%%%%
|[u(\cdot,t)>\mu_-+2\xi\om]\cap B_\rho|
\ge{\txty\frac14}|B_\rho|\quad\text{ for all }\> 
t \in(-\xi\om\rho,0].
%%%%%%%%%%%%%%%%%%%%%%%%%%%%%%%%%%%%%%%%%%%%%%%%%%
\end{equation}
%%%%%%%%%%%%%%%%%%%%%%%%%%%%%%%%%%%%%%%%%%%%%%%%%%
Next, apply the discrete isoperimetric inequality 
of Lemma~{2.2} of \cite[Preliminaries]{DBGV-mono} to 
the function $u(\cdot,t)$, for $t$ in the range 
$(-\xi\om\rho,0]$, over the ball $B_\rho$, for the levels 
%%%%%%%%%%%%%%%%%%%%%%%%%%%%%%%%%%%%%%%%%%%%%%%%%%
\begin{equation*}
%%%%%%%%%%%%%%%%%%%%%%%%%%%%%%%%%%%%%%%%%%%%%%%%%%
k=\mu_-+\xi\om \ \ \text{ and }\ \ 
\ell=\mu_-+2\xi\om\quad\text{ so that }\quad \ell-k
=\xi\om.
%%%%%%%%%%%%%%%%%%%%%%%%%%%%%%%%%%%%%%%%%%%%%%%%%%
\end{equation*}
%%%%%%%%%%%%%%%%%%%%%%%%%%%%%%%%%%%%%%%%%%%%%%%%%%
This inequality is stated and proved in \cite{DBGV-mono} 
for functions in $W^{1,1}_{\loc}(E)$. It continues 
to hold for $u\in BV_{\loc}(E)$, by virtue of 
the approximation procedure of \cite[Theorem~1.17]{giusti}.  
Taking also into account (\ref{Eq:8:2}) this gives
%%%%%%%%%%%%%%%%%%%%%%%%%%%%%%%%%%%%%%%%%%%%%%%%%%
\begin{equation*}
%%%%%%%%%%%%%%%%%%%%%%%%%%%%%%%%%%%%%%%%%%%%%%%%%%
\xi\om|[u(\cdot,t)<\mu_-+\xi\om]\cap B_\rho|
\le\gm\rho\|Du\|([u(\cdot,t)>k]\cap B_\rho).
%%%%%%%%%%%%%%%%%%%%%%%%%%%%%%%%%%%%%%%%%%%%%%%%%%
\end{equation*}
%%%%%%%%%%%%%%%%%%%%%%%%%%%%%%%%%%%%%%%%%%%%%%%%%%
Integrating in $dt$ over the time interval 
$(-\xi\om\rho,0]$, gives
%%%%%%%%%%%%%%%%%%%%%%%%%%%%%%%%%%%%%%%%%%%%%%%%%%
\begin{equation*}
%%%%%%%%%%%%%%%%%%%%%%%%%%%%%%%%%%%%%%%%%%%%%%%%%%
\frac{\Big|\big[u<\mu_-+\xi\om\big]\cap 
Q_\rho(\xi\om)\big|}{\big|Q_\rho(\xi\om)\big|} 
\le\frac{\gm}{(\xi\om_o)^2}\,\frac\rho{|Q_\rho|}\, 
\int_{-\rho}^0\|Du(\cdot,t)\|(B_\rho)dt. 
%%%%%%%%%%%%%%%%%%%%%%%%%%%%%%%%%%%%%%%%%%%%%%%%%%
\end{equation*}
%%%%%%%%%%%%%%%%%%%%%%%%%%%%%%%%%%%%%%%%%%%%%%%%%%
By the assumption, the right-hand side tends to zero 
as $\rho\searrow0$. 
%%%%%%%%%%%%%%%%%%%%%%%%%%%%%%%%%%%%%%%%%%%%%%%%%%
Hence, there exists $\rho$ so small that
%%%%%%%%%%%%%%%%%%%%%%%%%%%%%%%%%%%%%%%%%%%%%%%%%%
\begin{equation*}
%%%%%%%%%%%%%%%%%%%%%%%%%%%%%%%%%%%%%%%%%%%%%%%%%%
\frac{\Big|\big[u<\mu_-+\xi\om\big]\cap 
Q_{\rho}(\xi\om)\big|}{\big|Q_{\rho}(\xi\om)\big|} 
\le\nu_- 
%%%%%%%%%%%%%%%%%%%%%%%%%%%%%%%%%%%%%%%%%%%%%%%%%%
\end{equation*}
%%%%%%%%%%%%%%%%%%%%%%%%%%%%%%%%%%%%%%%%%%%%%%%%%%
where $\nu_-$ is the number claimed by Lemma~\ref{Lm:6:1} 
for such choice of parameters. The Lemma then implies 
%%%%%%%%%%%%%%%%%%%%%%%%%%%%%%%%%%%%%%%%%%%%%%%%%%
\begin{equation*}
%%%%%%%%%%%%%%%%%%%%%%%%%%%%%%%%%%%%%%%%%%%%%%%%%%
\essinf_{Q_{\frac12\rho}(\xi\om)}u\ge 
\mu_-+{\txty\frac12}\xi\om,
%%%%%%%%%%%%%%%%%%%%%%%%%%%%%%%%%%%%%%%%%%%%%%%%%%
\end{equation*}
%%%%%%%%%%%%%%%%%%%%%%%%%%%%%%%%%%%%%%%%%%%%%%%%%%
and hence
%%%%%%%%%%%%%%%%%%%%%%%%%%%%%%%%%%%%%%%%%%%%%%%%%%
\begin{equation*}
%%%%%%%%%%%%%%%%%%%%%%%%%%%%%%%%%%%%%%%%%%%%%%%%%%
\essosc_{Q_{\frac12\rho}(\xi\om)} u\le\eta\om
\qquad\text{ where }\quad 
\eta=1-{\txty\frac12}\xi\in (0,1).
%%%%%%%%%%%%%%%%%%%%%%%%%%%%%%%%%%%%%%%%%%%%%%%%%%
\end{equation*}
%%%%%%%%%%%%%%%%%%%%%%%%%%%%%%%%%%%%%%%%%%%%%%%%%%
Setting $\rho_1=\frac12\xi\om\rho$ gives 
%%%%%%%%%%%%%%%%%%%%%%%%%%%%%%%%%%%%%%%%%%%%%%%%%%
\begin{equation*}
%%%%%%%%%%%%%%%%%%%%%%%%%%%%%%%%%%%%%%%%%%%%%%%%%%
\om_{\rho_1}=\essosc_{Q_{\rho_1}} u\le\eta\om.
%%%%%%%%%%%%%%%%%%%%%%%%%%%%%%%%%%%%%%%%%%%%%%%%%%
\end{equation*}
%%%%%%%%%%%%%%%%%%%%%%%%%%%%%%%%%%%%%%%%%%%%%%%%%%
Repeat now the same argument starting from the 
cylinder $Q_{\rho_1}$, and proceed recursively to 
generate a decreasing sequence of radii  
$\{\rho_n\}\to0$ such that 
%%%%%%%%%%%%%%%%%%%%%%%%%%%%%%%%%%%%%%%%%%%%%%%%%%
\begin{equation*}
%%%%%%%%%%%%%%%%%%%%%%%%%%%%%%%%%%%%%%%%%%%%%%%%%%
\om_o\le \essosc_{Q_{\rho_n}} u\le\eta^n\om
\quad\text{ for all }\> n\in\nn.\tag*{\bbox}
%%%%%%%%%%%%%%%%%%%%%%%%%%%%%%%%%%%%%%%%%%%%%%%%%%
\end{equation*}
%%%%%%%%%%%%%%%%%%%%%%%%%%%%%%%%%%%%%%%%%%%%%%%%%%
\AppendicesFromNowOn
%%%%%%%%%%%%%%%%%%%%%%%%%%%%%%%%%%%%%%%%%%%%%%%%%%
\Appendix{Proof of Proposition~\ref{Prop:3:1}}\label{App:A}
%%%%%%%%%%%%%%%%%%%%%%%%%%%%%%%%%%%%%%%%%%%%%%%%%%
The proof uses an approximation procedure of \cite{BDKM}. 
Observe first that the assumption $u_t\in L^1_{\loc}(E_T)$ 
permits to cast (\ref{Eq:1:1}) in the form
%%%%%%%%%%%%%%%%%%%%%%%%%%%%%%%%%%%%%%%%%%%%%%%%%%%
\begin{equation}\label{Eq:A:1}
%%%%%%%%%%%%%%%%%%%%%%%%%%%%%%%%%%%%%%%%%%%%%%%%%%%
\|Du(t)\|(E)\le\|D(u+\vp)(t)\|(E)-\ine u_t\vp dx
%%%%%%%%%%%%%%%%%%%%%%%%%%%%%%%%%%%%%%%%%%%%%%%%%%%
\end{equation}
%%%%%%%%%%%%%%%%%%%%%%%%%%%%%%%%%%%%%%%%%%%%%%%%%%%
for a.e. $t\in(0,T)$ for all 
%%%%%%%%%%%%%%%%%%%%%%%%%%%%%%%%%%%%%%%%%%%%%%%%%%%
\begin{equation}\label{Eq:A:2}
%%%%%%%%%%%%%%%%%%%%%%%%%%%%%%%%%%%%%%%%%%%%%%%%%%%
\vp\in BV_{\loc}(E)\cap L^\infty_{\loc}(E)
\quad\text{ with }\> \supp\{\vp\}\subset E. 
%%%%%%%%%%%%%%%%%%%%%%%%%%%%%%%%%%%%%%%%%%%%%%%%%%%
\end{equation}
%%%%%%%%%%%%%%%%%%%%%%%%%%%%%%%%%%%%%%%%%%%%%%%%%%%
We only prove the estimate for $\ukp$, 
the one for $\ukm$ being similar. Fix a cylinder 
%%%%%%%%%%%%%%%%%%%%%%%%%%%%%%%%%%%%%%%%%%%%%%%%%%%
\begin{equation*}
%%%%%%%%%%%%%%%%%%%%%%%%%%%%%%%%%%%%%%%%%%%%%%%%%%%
\big[\pto+Q_{\rho}(\theta)\big]\subset E_T. 
%%%%%%%%%%%%%%%%%%%%%%%%%%%%%%%%%%%%%%%%%%%%%%%%%%%
\end{equation*}
%%%%%%%%%%%%%%%%%%%%%%%%%%%%%%%%%%%%%%%%%%%%%%%%%%%
Up to a translation, assume that $\pto=(0,0)$ and 
fix a time $t\in (-\theta\rho,0)$ for which 
%%%%%%%%%%%%%%%%%%%%%%%%%%%%%%%%%%%%%%%%%%%%%%%%%%%
\begin{equation*}
%%%%%%%%%%%%%%%%%%%%%%%%%%%%%%%%%%%%%%%%%%%%%%%%%%%
\int_{B_{\rho}}|u_t(x,t)|dx<\infty, 
\quad\text{ and }\quad 
u(\cdot, t)\in BV(E)\cap L^\infty(B_\rho).
%%%%%%%%%%%%%%%%%%%%%%%%%%%%%%%%%%%%%%%%%%%%%%%%%%%
\end{equation*}
%%%%%%%%%%%%%%%%%%%%%%%%%%%%%%%%%%%%%%%%%%%%%%%%%%%
The next approximation procedure is carried {out} for 
such $t$ fixed and we write $u(\cdot,t)=u$. 
By \cite[Theorem~1.17]{giusti}, there exists 
$\{u_j\}\subset C^\infty(B_{\rho})$ such that
%%%%%%%%%%%%%%%%%%%%%%%%%%%%%%%%%%%%%%%%%%%%%%%%%%
\begin{equation}\label{Eq:A:3}
%%%%%%%%%%%%%%%%%%%%%%%%%%%%%%%%%%%%%%%%%%%%%%%%%%
\lim_{j\to\infty}\int_{B_{\rho}}|u_j-u|dx=0
\quad\text{ and }\quad 
\|Du\|(E)=\lim_{j\to\infty}\ine|Du_j|dx. 
%%%%%%%%%%%%%%%%%%%%%%%%%%%%%%%%%%%%%%%%%%%%%%%%%%
\end{equation} 
%%%%%%%%%%%%%%%%%%%%%%%%%%%%%%%%%%%%%%%%%%%%%%%%%%
Test (\ref{Eq:A:1}) with $\vp=-\z(u-k)_+$, 
where {$\z\in\mcl{C}\big(Q_\rho(\theta)\big)$}. This is an admissible choice, 
since $u\in BV(E)\cap L^\infty(B_\rho)$. Set 
$\vp_j=-\z(u_j-k)_+$ for $j\in\nn$. For a given $\eps>0$ 
there exists $j_o\in\nn$ such that
%%%%%%%%%%%%%%%%%%%%%%%%%%%%%%%%%%%%%%%%%%%%%%%%%%
\begin{equation*}%%\label{Eq:A:3}
%%%%%%%%%%%%%%%%%%%%%%%%%%%%%%%%%%%%%%%%%%%%%%%%%%
\int_{E}|Du_j|dx<\|Du(\cdot,t)\|(E)+\frac12\eps
\quad\text{ for all }\>j\ge j_o.
%%%%%%%%%%%%%%%%%%%%%%%%%%%%%%%%%%%%%%%%%%%%%%%%%%
\end{equation*} 
%%%%%%%%%%%%%%%%%%%%%%%%%%%%%%%%%%%%%%%%%%%%%%%%%%
Here we have used the second of (\ref{Eq:A:3}). 
By the first, $\{(u_j+\vp_j)\}\to (u+\vp)$  
in $L^1(E)$. Therefore, for any $\bpsi \in [C_o^1(E)]^N$ 
with $\|\bpsi\|\le1$, 
%%%%%%%%%%%%%%%%%%%%%%%%%%%%%%%%%%%%%%%%%%%%%%%%%%
\begin{align*}
%%%%%%%%%%%%%%%%%%%%%%%%%%%%%%%%%%%%%%%%%%%%%%%%%%
\ine(u+\vp)\dvg\bpsi\, dx&=\lim_{j\to\infty}\ine(u_j+\vp_j)\dvg\bpsi\, dx\\
&\le\liminf_{j\to\infty}\ine|D(u_j +\vp_j)|\,dx.
%%%%%%%%%%%%%%%%%%%%%%%%%%%%%%%%%%%%%%%%%%%%%%%%%%
\end{align*}
%%%%%%%%%%%%%%%%%%%%%%%%%%%%%%%%%%%%%%%%%%%%%%%%%%
Taking the supremum over all such $\bpsi$ gives
%%%%%%%%%%%%%%%%%%%%%%%%%%%%%%%%%%%%%%%%%%%%%%%%%%
\begin{equation*}%%\label{Eq:A:3}
%%%%%%%%%%%%%%%%%%%%%%%%%%%%%%%%%%%%%%%%%%%%%%%%%%
\|D(u+\vp)(t)\|(E)\le\liminf_{j\to\infty}\ine|D(u_j +\vp_j)|\,dx.
%%%%%%%%%%%%%%%%%%%%%%%%%%%%%%%%%%%%%%%%%%%%%%%%%%
\end{equation*} 
%%%%%%%%%%%%%%%%%%%%%%%%%%%%%%%%%%%%%%%%%%%%%%%%%%
Therefore, up to redefining $j_o$ we may also assume that
%%%%%%%%%%%%%%%%%%%%%%%%%%%%%%%%%%%%%%%%%%%%%%%%%%
\begin{equation*}%%\label{Eq:A:3}
%%%%%%%%%%%%%%%%%%%%%%%%%%%%%%%%%%%%%%%%%%%%%%%%%%
\ine|D(u_j+\vp_j)|dx\ge\|D(u+\vp )\|(E)-\frac12\eps
\quad\text{for all}\ \ j\ge j_o.
%%%%%%%%%%%%%%%%%%%%%%%%%%%%%%%%%%%%%%%%%%%%%%%%%%
\end{equation*} 
%%%%%%%%%%%%%%%%%%%%%%%%%%%%%%%%%%%%%%%%%%%%%%%%%%
Combining the preceding inequalities gives that
%%%%%%%%%%%%%%%%%%%%%%%%%%%%%%%%%%%%%%%%%%%%%%%%%%
\begin{align}\label{Eq:A:4}
%%%%%%%%%%%%%%%%%%%%%%%%%%%%%%%%%%%%%%%%%%%%%%%%%%
\ine|Du_j|dx&<\|Du(\cdot,t)\|(E)+\frac12\eps\notag\\
&\le\|D(u+\vp)(\cdot,t)\|(E)+\ine u_t(\cdot,t)\vp dx+\frac12\eps\\
&\le\ine|D(u_j+\vp_j)|dx+\ine u_t(\cdot,t)\vp dx+\eps\notag
%%%%%%%%%%%%%%%%%%%%%%%%%%%%%%%%%%%%%%%%%%%%%%%%%%
\end{align}
%%%%%%%%%%%%%%%%%%%%%%%%%%%%%%%%%%%%%%%%%%%%%%%%%%
for all $j\ge j_o$. Next, estimate the first integral 
on the right-hand side as,
%%%%%%%%%%%%%%%%%%%%%%%%%%%%%%%%%%%%%%%%%%%%%%%%%%
\begin{align*}
%%%%%%%%%%%%%%%%%%%%%%%%%%%%%%%%%%%%%%%%%%%%%%%%%%
\ine&|D(u_j +\vp_j)|dx=\ine|D(u_j-\z(u_j-k)_+)|dx\\
&\le\ine|Du_j-\z D(u_j-k)_+|dx+\ine|D\z|(u_j-k)_+dx\\
&\le \ine(1-\z)|Du_j|+\z|Du_j-D(u_j-k)_+|dx
+\ine|D\z|(u_j-k)_+dx.
%%%%%%%%%%%%%%%%%%%%%%%%%%%%%%%%%%%%%%%%%%%%%%%%%%
\end{align*}
%%%%%%%%%%%%%%%%%%%%%%%%%%%%%%%%%%%%%%%%%%%%%%%%%%
Put this in (\ref{Eq:A:4}), and absorb the first integral 
on the right-hand side into the left-hand side, to obtain 
%%%%%%%%%%%%%%%%%%%%%%%%%%%%%%%%%%%%%%%%%%%%%%%%%%
\begin{align*}
%%%%%%%%%%%%%%%%%%%%%%%%%%%%%%%%%%%%%%%%%%%%%%%%%%
\ine\z|D(u_j-k)_+|dx&=\ine\z\big[|Du_j|-|Du_j
-D(u_j-k)_+|\big]dx\\
&\le\ine|D\z|(u_j-k)_+dx+\ine u_t(\cdot,t)\vp dx+\eps.
%%%%%%%%%%%%%%%%%%%%%%%%%%%%%%%%%%%%%%%%%%%%%%%%%%
\end{align*}
%%%%%%%%%%%%%%%%%%%%%%%%%%%%%%%%%%%%%%%%%%%%%%%%%%
From this 
%%%%%%%%%%%%%%%%%%%%%%%%%%%%%%%%%%%%%%%%%%%%%%%%%%
\begin{equation*}
%%%%%%%%%%%%%%%%%%%%%%%%%%%%%%%%%%%%%%%%%%%%%%%%%%
\ine|D(\z(u_j-k)_+)|dx\le2\ine|D\z|(u_j-k)_+dx
+\ine u_t(\cdot,t)\vp dx+\eps.
%%%%%%%%%%%%%%%%%%%%%%%%%%%%%%%%%%%%%%%%%%%%%%%%%%
\end{equation*}
%%%%%%%%%%%%%%%%%%%%%%%%%%%%%%%%%%%%%%%%%%%%%%%%%%
Next let $j\to\infty$, using the lower semicontinuity 
of the total variation with respect to $L^1$-convergence. 
This gives 
%%%%%%%%%%%%%%%%%%%%%%%%%%%%%%%%%%%%%%%%%%%%%%%%%%
\begin{align*}
%%%%%%%%%%%%%%%%%%%%%%%%%%%%%%%%%%%%%%%%%%%%%%%%%%
\|D(\z(u-k)_+)\|(B_\rho)
&\le\liminf_{j\to\infty}\ine|D(\z(u_j-k)_+)|dx\\
&\le\lim_{j\to\infty}2\ine|D\z|(u_j-k)_+dx+\ine u_t\vp dx+\eps\\
&=2\ine|D\z|(u-k)_+dx+\ine u_t\vp dx+\eps.
%%%%%%%%%%%%%%%%%%%%%%%%%%%%%%%%%%%%%%%%%%%%%%%%%%
\end{align*}
%%%%%%%%%%%%%%%%%%%%%%%%%%%%%%%%%%%%%%%%%%%%%%%%%%
Finally let $\eps\to 0$ and use the definition of $\vp$ to get
%%%%%%%%%%%%%%%%%%%%%%%%%%%%%%%%%%%%%%%%%%%%%%%%%%
\begin{equation*}
%%%%%%%%%%%%%%%%%%%%%%%%%%%%%%%%%%%%%%%%%%%%%%%%%%
\|D(\z(u-k)_+)\|(B_\rho)\le2\int_{B_\rho}|D\z|(u-k)_+dx-\int_{B_\rho}\z u_t (u-k)_+dx.
%%%%%%%%%%%%%%%%%%%%%%%%%%%%%%%%%%%%%%%%%%%%%%%%%%
\end{equation*}
%%%%%%%%%%%%%%%%%%%%%%%%%%%%%%%%%%%%%%%%%%%%%%%%%%
To conclude the proof, integrate in $dt$ over $(-\theta\rho,0)$.\hfill\bbox
%%%%%%%%%%%%%%%%%%%%%%%%%%%%%%%%%%%%%%%%%%%%%%%%%%%%%
\vskip.2truecm
%%%%%%%%%%%%%%%%%%%%%%%%%%%%%%%%%%%%%%%%%%%%%%%%%%%%%
\Appendix{Boundedness of Minimizers}\label{App:B}
%%%%%%%%%%%%%%%%%%%%%%%%%%%%%%%%%%%%%%%%%%%%%%%%%%%%%
%%%%%%%%%%%%%%%%%%%%%%%%%%%%%%%%%%%%%%%%%%%%%%%%%%%%%
\begin{proposition}\label{Prop:B:1} 
%%%%%%%%%%%%%%%%%%%%%%%%%%%%%%%%%%%%%%%%%%%%%%%%%%%%%
Let $u:E_T\to\rr$ be a parabolic minimizer of the total 
variation flow in the sense of (\ref{Eq:1:1}). Furthermore, 
assume that $u\in L^r_{\loc}(E_T)$ for some  $r>N$, and that 
it can be constructed as the limit in $L^r_{\loc}(E_T)$ 
of a sequence of parabolic minimizers satisfying 
(\ref{Eq:1:2}).  Then, there exists a positive 
constant $\gm$ depending only upon $N,\gm,r$, such that
%%%%%%%%%%%%%%%%%%%%%%%%%%%%%%%%%%%%%%%%%%%%%%%%%%%%%
\begin{equation}\label{Eq:B:1}
%%%%%%%%%%%%%%%%%%%%%%%%%%%%%%%%%%%%%%%%%%%%%%%%%%%%%
\begin{aligned}
%%%%%%%%%%%%%%%%%%%%%%%%%%%%%%%%%%%%%%%%%%%%%%%%%%%%%
\sup_{B_\rho(y)\times[s,t]}u_\pm&\le
\gm\Big(\frac\rho{t-s}\Big)^{\frac{N}{r-N}}
\Big(\frac1{\rho^N(t-s)}\int_{2s-t}^t\int_{B_{4\rho}(y)}u_\pm^r\,dxd\tau
\Big)^{\frac{1}{r-N}}\\
&+\gm\frac{t-s}{\rho}
%%%%%%%%%%%%%%%%%%%%%%%%%%%%%%%%%%%%%%%%%%%%%%%%%%%%%
\end{aligned}
%%%%%%%%%%%%%%%%%%%%%%%%%%%%%%%%%%%%%%%%%%%%%%%%%%%%%
\end{equation}
%%%%%%%%%%%%%%%%%%%%%%%%%%%%%%%%%%%%%%%%%%%%%%%%%%%%%
for all cylinders
%%%%%%%%%%%%%%%%%%%%%%%%%%%%%%%%%%%%%%%%%%%%%%%%%%%%%
\begin{equation*}
%%%%%%%%%%%%%%%%%%%%%%%%%%%%%%%%%%%%%%%%%%%%%%%%%%%%%
B_{4\rho}(y)\times[s-(t-s),s+(t-s)]\subset E_T.
%%%%%%%%%%%%%%%%%%%%%%%%%%%%%%%%%%%%%%%%%%%%%%%%%%%%%
\end{equation*}
%%%%%%%%%%%%%%%%%%%%%%%%%%%%%%%%%%%%%%%%%%%%%%%%%%%%%
The constant $\gm(N,\gm,r)\to\infty$ as either 
$r\to N$, or $r\to\infty$.
%%%%%%%%%%%%%%%%%%%%%%%%%%%%%%%%%%%%%%%%%%%%%%%%%%%%%
\end{proposition}
%%%%%%%%%%%%%%%%%%%%%%%%%%%%%%%%%%%%%%%%%%%%%%%%%%%%%%%%%%%%%
\begin{remark}\label{Rmk:B:1} {\normalfont
%%%%%%%%%%%%%%%%%%%%%%%%%%%%%%%%%%%%%%%%%%%%%%%%%%%%%%%%%%%%%
It is not required that the approximations to $u$ satisfy 
(\ref{Eq:1:2}) uniformly. The latter is only needed to cast 
a function satisfying (\ref{Eq:1:1}) into a DeGiorgi class. 
The proof of the proposition only uses such a membership, and 
turns such a {\it qualitative}, non-uniform information 
into the {\it quantitative} information (\ref{Eq:B:1}).
%%%%%%%%%%%%%%%%%%%%%%%%%%%%%%%%%%%%%%%%%%%%%%%%%%%%%%%%%%%%%
}%%
%%%%%%%%%%%%%%%%%%%%%%%%%%%%%%%%%%%%%%%%%%%%%%%%%%%%%%%%%%%%%
\end{remark}
%%%%%%%%%%%%%%%%%%%%%%%%%%%%%%%%%%%%%%%%%%%%%%%%%%%%%%%%%%%%%
\noi{\bf Proof} (of Proposition~\ref{Prop:B:1}). Let $\{u_j\}$ 
be a sequence of approximating functions to $u$. Since 
$u_j$ satisfy (\ref{Eq:1:2}), they belong to the 
classes $[DG](E_T;2)$, by Proposition~\ref{Prop:3:1}. It will 
suffice to establish (\ref{Eq:B:1}) for such $u_j$ for a 
constant $\gm$ independent of $j$. Thus in the calculations below 
we drop the suffix {$j$} from $u_j$. The proof will be given 
for non-negative $u\in [DG]^+(E_T;2)$, the proof for the 
remaining case being identical; {it is very similar to the proof of 
Proposition~A.2.1 given in \cite[\S~A.2]{DBGV-mono}}.
Assume $(y,s)=(0,0)$ and for fixed 
$\sig\in(0,1)$ and $n=0,1,2,\dots$ set
%%%%%%%%%%%%%%%%%%%%%%%%%%%%%%%%%%%%%%%%%%%%%%%%%%%%%%%%%%%%%
\begin{equation*}
%%%%%%%%%%%%%%%%%%%%%%%%%%%%%%%%%%%%%%%%%%%%%%%%%%%%%%%%%%%%%
\begin{array}{ll}
%%%%%%%%%%%%%%%%%%%%%%%%%%%%%%%%%%%%%%%%%%%%%%%%%%%%%%%%%%%%%
{\dsty \rho_n=\sig\rho+\frac{1-\sig}{2^n}\rho,\qquad }
&{\dsty t_n=-\sig t-\frac{1-\sig}{2^n}t},\\
{\dsty B_n=B_{\rho_n},\qquad }&{\dsty Q_n=B_n\times(t_n,t).}
%%%%%%%%%%%%%%%%%%%%%%%%%%%%%%%%%%%%%%%%%%%%%%%%%%%%%%%%%%%%%
\end{array}
%%%%%%%%%%%%%%%%%%%%%%%%%%%%%%%%%%%%%%%%%%%%%%%%%%%%%%%%%%%%%
\end{equation*}
%%%%%%%%%%%%%%%%%%%%%%%%%%%%%%%%%%%%%%%%%%%%%%%%%%%%%%%%%%%%%
This is a family of nested and shrinking cylinders with common 
``vertex'' at $(0,t)$, and by construction 
%%%%%%%%%%%%%%%%%%%%%%%%%%%%%%%%%%%%%%%%%%%%%%%%%%%%%%%%%%%%%
\begin{equation*}
%%%%%%%%%%%%%%%%%%%%%%%%%%%%%%%%%%%%%%%%%%%%%%%%%%%%%%%%%%%%%
Q_o=B_\rho\times(-t,t)\quad\text{ and }\quad 
Q_\infty=B_{\sig\rho}\times(-\sig t,t).
%%%%%%%%%%%%%%%%%%%%%%%%%%%%%%%%%%%%%%%%%%%%%%%%%%%%%%%%%%%%%
\end{equation*}
%%%%%%%%%%%%%%%%%%%%%%%%%%%%%%%%%%%%%%%%%%%%%%%%%%%%%%%%%%%%%
We have assumed that $u$ can be constructed as the 
limit in $L^r_{\loc}(E_T)$ of a sequence of bounded 
parabolic minimizers. By working with such approximations, 
we may assume that $u$ is qualitatively locally bounded. 
Therefore, set 
%%%%%%%%%%%%%%%%%%%%%%%%%%%%%%%%%%%%%%%%%%%%%%%%%%%%%%%%%%%%%
\begin{equation*}
%%%%%%%%%%%%%%%%%%%%%%%%%%%%%%%%%%%%%%%%%%%%%%%%%%%%%%%%%%%%%
M=\essup_{Q_o}\max\{u,0\},\qquad M_\sig
=\essup_{Q_\infty}\max\{u,0\}.
%%%%%%%%%%%%%%%%%%%%%%%%%%%%%%%%%%%%%%%%%%%%%%%%%%%%%%%%%%%%%
\end{equation*}
%%%%%%%%%%%%%%%%%%%%%%%%%%%%%%%%%%%%%%%%%%%%%%%%%%%%%%%%%%%%%
We first find a relationship between $M$ and $M_\sig$. 
Denote by $\z$ a non-negative, piecewise smooth cutoff 
function in $Q_n$ that equals one on $Q_{n+1}$, and has the form 
$\z(x,t)=\z_1(x)\z_2(t)$, where
%%%%%%%%%%%%%%%%%%%%%%%%%%%%%%%%%%%%%%%%%%%%%%%%%%%%%%%%%%%%%
\begin{equation*}
%%%%%%%%%%%%%%%%%%%%%%%%%%%%%%%%%%%%%%%%%%%%%%%%%%%%%%%%%%%%%
\begin{array}{lc}
%%%%%%%%%%%%%%%%%%%%%%%%%%%%%%%%%%%%%%%%%%%%%%%%%%%%%%%%%%%%%
{\dsty 
%%%%%%%%%%%%%%%%%%%%%%%%%%%%%%%%%%%%%%%%%%%%%%%%%%%%%%%%%%%%%
\z_1=\left\{
%%%%%%%%%%%%%%%%%%%%%%%%%%%%%%%%%%%%%%%%%%%%%%%%%%%%%%%%%%%%%
\begin{array}{ll}
%%%%%%%%%%%%%%%%%%%%%%%%%%%%%%%%%%%%%%%%%%%%%%%%%%%%%%%%%%%%%
1\>&\text{ in }\> B_{n+1}\\
0\>&\text{ in }\>\rn-B_n
%%%%%%%%%%%%%%%%%%%%%%%%%%%%%%%%%%%%%%%%%%%%%%%%%%%%%%%%%%%%%
\end{array}\right .}\>
%%%%%%%%%%%%%%%%%%%%%%%%%%%%%%%%%%%%%%%%%%%%%%%%%%%%%%%%%%%%%
&{\dsty |D\z_1|\le\frac{2^{n+1}}{(1-\sig)\rho}}\\
%%%%%%%%%%%%%%%%%%%%%%%%%%%%%%%%%%%%%%%%%%%%%%%%%%%%%%%%%%%%%
{}\\
%%%%%%%%%%%%%%%%%%%%%%%%%%%%%%%%%%%%%%%%%%%%%%%%%%%%%%%%%%%%%
{\dsty 
%%%%%%%%%%%%%%%%%%%%%%%%%%%%%%%%%%%%%%%%%%%%%%%%%%%%%%%%%%%%%
\z_2=\left\{
%%%%%%%%%%%%%%%%%%%%%%%%%%%%%%%%%%%%%%%%%%%%%%%%%%%%%%%%%%%%%
\begin{array}{ll}
%%%%%%%%%%%%%%%%%%%%%%%%%%%%%%%%%%%%%%%%%%%%%%%%%%%%%%%%%%%%%
0\>&\text{ for }\> t\le t_n\\
1\>&\text{ for }\> t\ge t_{n+1}
%%%%%%%%%%%%%%%%%%%%%%%%%%%%%%%%%%%%%%%%%%%%%%%%%%%%%%%%%%%%%
\end{array}\right .}\>
%%%%%%%%%%%%%%%%%%%%%%%%%%%%%%%%%%%%%%%%%%%%%%%%%%%%%%%%%%%%%
&{\dsty 0\le \z_{2,t}\le\frac{2^{n+1}}{(1-\sig)t}};
%%%%%%%%%%%%%%%%%%%%%%%%%%%%%%%%%%%%%%%%%%%%%%%%%%%%%%%%%%%%%
\end{array}
%%%%%%%%%%%%%%%%%%%%%%%%%%%%%%%%%%%%%%%%%%%%%%%%%%%%%%%%%%%%%
\end{equation*} 
%%%%%%%%%%%%%%%%%%%%%%%%%%%%%%%%%%%%%%%%%%%%%%%%%%%%%%%%%%%%%
introduce the increasing sequence of levels $k_n=k-2^{-n}k$, 
where $k>0$ is to be  chosen, and in (\ref{Eq:3:2})$_+$, 
take such a test function, to get
%%%%%%%%%%%%%%%%%%%%%%%%%%%%%%%%%%%%%%%%%%%%%%%%%%%%%%%%%%%%%
\begin{align}\label{Eq:B:2}
%%%%%%%%%%%%%%%%%%%%%%%%%%%%%%%%%%%%%%%%%%%%%%%%%%%%%%%%%%%%%
\sup_{t_n\le \tau\le t}&\int_{B_n}[\ukno\z]^2(x,\tau)dx+
\int_{t_n}^t\|D [\ukno\z](\cdot,\tau)\|(B_n)\,d\tau\nonumber\\
&\le\frac{\gm 2^{n}}{(1-\sig)\rho}\iint_{Q_n}\uknpu\,dx\,d\tau\\
&+\frac{\gm 2^{n}}{(1-\sig)t}\iint_{Q_n}\uknpu^2dx\,d\tau.\nonumber
%%%%%%%%%%%%%%%%%%%%%%%%%%%%%%%%%%%%%%%%%%%%%%%%%%%%%%%%%%%%%
\end{align}
%%%%%%%%%%%%%%%%%%%%%%%%%%%%%%%%%%%%%%%%%%%%%%%%%%%%%%%%%%%%%%%%%%%%
%The requirement $\lm_r>0$ implies 
%\begin{equation*}
%r>2\ge q=p{\txty\frac{N+2}{N}}.
%\end{equation*}
Estimate 
%%%%%%%%%%%%%%%%%%%%%%%%%%%%%%%%%%%%%%%%%%%%%%%%%%%%%%%%%%%%%
\begin{align*}
%%%%%%%%%%%%%%%%%%%%%%%%%%%%%%%%%%%%%%%%%%%%%%%%%%%%%%%%%%%%%
\iint_{Q_n}\uknpu dx\,d\tau&\le\gm\frac{2^{n(r-1)}}{k^{r-1}}
\iint_{Q_n}\uknp^rdx\,d\tau,\\
\iint_{Q_n}\uknpu^2dx\,d\tau&\le\gm\frac{2^{n(r-2)}}{k^{r-2}}
\iint_{Q_n}\uknp^rdx\,d\tau.\\
%%%%%%%%%%%%%%%%%%%%%%%%%%%%%%%%%%%%%%%%%%%%%%%%%%%%%%%%%%%%%
\end{align*}
%%%%%%%%%%%%%%%%%%%%%%%%%%%%%%%%%%%%%%%%%%%%%%%%%%%%%%%%%%%%%
Taking these estimates into account yields
%%%%%%%%%%%%%%%%%%%%%%%%%%%%%%%%%%%%%%%%%%%%%%%%%%%%%%%%%%%%%
\begin{align*}
%%%%%%%%%%%%%%%%%%%%%%%%%%%%%%%%%%%%%%%%%%%%%%%%%%%%%%%%%%%%%
&\sup_{t_n<\tau\le t}\int_{B_{n}}[\uknpu\zeta]^2(x,\tau)dx
+\int_{t_n}^t\|D[\uknpu\zeta](\cdot,\tau)\|(B_n)\,d\tau\\
&\le\gm\frac{2^{nr}}{(1-\sig) t}\Big[
\Big(\frac{t}{\rho}\Big)k^{1-r}+\frac{1}{k^{r-2}}\Big]
\iint_{Q_n}\uknp^rdx\,d\tau.
%%%%%%%%%%%%%%%%%%%%%%%%%%%%%%%%%%%%%%%%%%%%%%%%%%%%%%%%%%%%%
\end{align*}
%%%%%%%%%%%%%%%%%%%%%%%%%%%%%%%%%%%%%%%%%%%%%%%%%%%%%%%%%%%%%
Assuming that $\dsty k>\frac t\rho$, this implies
%%%%%%%%%%%%%%%%%%%%%%%%%%%%%%%%%%%%%%%%%%%%%%%%%%%%%%%%%%%%%
\begin{align*}
%%%%%%%%%%%%%%%%%%%%%%%%%%%%%%%%%%%%%%%%%%%%%%%%%%%%%%%%%%%%%
\sup_{t_n<\tau\le t}\int_{B_{n}}[\uknpu\z]^2&(x,\tau)dx
+\int_{t_n}^t\|D[\uknpu\z](\cdot,\tau)\|(B_n)\,d\tau\\
&\le\frac{\gm 2^{nr}}{(1-\sig) t}\frac{1}{k^{r-2}}
\iint_{Q_n}\uknp^rdx\,d\tau.
%%%%%%%%%%%%%%%%%%%%%%%%%%%%%%%%%%%%%%%%%%%%%%%%%%%%%%%%%%%%%
\end{align*}
%%%%%%%%%%%%%%%%%%%%%%%%%%%%%%%%%%%%%%%%%%%%%%%%%%%%%%%%%%%%%
Set 
%%%%%%%%%%%%%%%%%%%%%%%%%%%%%%%%%%%%%%%%%%%%%%%%%%%%%%%%%%%%%
\begin{equation*}
%%%%%%%%%%%%%%%%%%%%%%%%%%%%%%%%%%%%%%%%%%%%%%%%%%%%%%%%%%%%%
Y_n=\frac{1}{|Q_n|}\iint_{Q_n}(u-k_n)_+^rdx\,d\tau
%%%%%%%%%%%%%%%%%%%%%%%%%%%%%%%%%%%%%%%%%%%%%%%%%%%%%%%%%%%%%
\end{equation*}
%%%%%%%%%%%%%%%%%%%%%%%%%%%%%%%%%%%%%%%%%%%%%%%%%%%%%%%%%%%%%
and estimate
%%%%%%%%%%%%%%%%%%%%%%%%%%%%%%%%%%%%%%%%%%%%%%%%%%%%%%%%%%%%%
\begin{equation*}
%%%%%%%%%%%%%%%%%%%%%%%%%%%%%%%%%%%%%%%%%%%%%%%%%%%%%%%%%%%%%
Y_{n+1}\le\|u\|_{\infty,Q_o}^{r-q}\Big(\frac{1}{|Q_n|}
\iint_{Q_n}\uknpu^q dx\,d\tau\Big),
%%%%%%%%%%%%%%%%%%%%%%%%%%%%%%%%%%%%%%%%%%%%%%%%%%%%%%%%%%%%%
\end{equation*}
%%%%%%%%%%%%%%%%%%%%%%%%%%%%%%%%%%%%%%%%%%%%%%%%%%%%%%%%%%%%%
where $q\df=\frac{N+2}N$. Applying the embedding 
Proposition~{4.1} of \cite[Preliminaries]{DBGV-mono}, 
the previous inequality can be rewritten as
%%%%%%%%%%%%%%%%%%%%%%%%%%%%%%%%%%%%%%%%%%%%%%%%%%%%%%%%%%%%%
\begin{equation*}
%%%%%%%%%%%%%%%%%%%%%%%%%%%%%%%%%%%%%%%%%%%%%%%%%%%%%%%%%%%%%
Y_{n+1}\le\gm\|u\|_{\infty,Q_o}^{r-q}
\Big(\frac{\rho}{t}\Big)\frac{b^n}{(1-\sig)^{\frac{1}{N}(N+1)}}
\frac{1}{k^{(r-2)\frac{N+1}{N}}}Y_n^{1+\frac{1}{N}},
%%%%%%%%%%%%%%%%%%%%%%%%%%%%%%%%%%%%%%%%%%%%%%%%%%%%%%%%%%%%%
\end{equation*}
%%%%%%%%%%%%%%%%%%%%%%%%%%%%%%%%%%%%%%%%%%%%%%%%%%%%%%%%%%%%%
where $b=2^{r\frac{N+1}{N}}$.  Apply Lemma~{5.1} of 
\cite[Preliminaries]{DBGV-mono}, and conclude that $Y_n\to0$ as $n\to+\infty$, 
provided $k$ is chosen to satisfy
%%%%%%%%%%%%%%%%%%%%%%%%%%%%%%%%%%%%%%%%%%%%%%%%%%%%%%%%%%%%%
\begin{equation*}
%%%%%%%%%%%%%%%%%%%%%%%%%%%%%%%%%%%%%%%%%%%%%%%%%%%%%%%%%%%%%
Y_o=\tvls{Q_o}{24}\kern0.1cm u^rdx\,d\tau=
\gm(1-\sig)^{N+1}\|u\|_{\infty,Q_o}^{-(r-q){N}}
\Big(\frac{t}{\rho}\Big)^{N}k^{(r-2)({N+1})},
%%%%%%%%%%%%%%%%%%%%%%%%%%%%%%%%%%%%%%%%%%%%%%%%%%%%%%%%%%%%%
\end{equation*}
%%%%%%%%%%%%%%%%%%%%%%%%%%%%%%%%%%%%%%%%%%%%%%%%%%%%%%%%%%%%%
which yields
%%%%%%%%%%%%%%%%%%%%%%%%%%%%%%%%%%%%%%%%%%%%%%%%%%%%%%%%%%%%%
\begin{equation*}
%%%%%%%%%%%%%%%%%%%%%%%%%%%%%%%%%%%%%%%%%%%%%%%%%%%%%%%%%%%%%
M_\sig\le\tilde\gm\frac{M^{\frac{N(r-q)}{(N+1)(r-2)}}}{(1-\sig)^{\frac{1}{r-2}}}
\Big(\frac{\rho}{t}\Big)^{\frac{N}{(N+1)(r-2)}}
\Big(\tvls{Q_o}{24}\kern0.1cm u^r\,dx\,d\tau
\Big)^{\frac{1}{(r-2)(N+1)}}.
%%%%%%%%%%%%%%%%%%%%%%%%%%%%%%%%%%%%%%%%%%%%%%%%%%%%%%%%%%%%%
\end{equation*}
%%%%%%%%%%%%%%%%%%%%%%%%%%%%%%%%%%%%%%%%%%%%%%%%%%%%%%%%%%%%%
The proof is concluded by the interpolation Lemma~{5.2} of 
\cite[Preliminaries]{DBGV-mono}.\hfill\bbox
%%%%%%%%%%%%%%%%%%%%%%%%%%%%%%%%%%%%%%%%%%%%%%%%%%%%%%%%%%%%%
%%%%%%%%%%%%%%%%%%%%%%%%%%%%%%%%%%%%%

%%%%%%%%%%%%%%%%%%%%%%%%%%%%%%%%%%%%%%%%%%%%%%%%%%%%%
\bye